\documentclass[11pt]{amsart}
\usepackage{latexsym}
\usepackage{amsmath}
\usepackage{amssymb}

\newcommand{\eproof}{\mbox{\ }\hfill $\Box$ \par \vskip 10pt}

\newtheorem{Theorem}{Theorem}[section]
\newtheorem{lemma}[Theorem]{Lemma}
\newtheorem{prop}[Theorem]{Proposition}
\newtheorem{rem}[Theorem]{Remark}

\numberwithin{equation}{section}

\def\cal{\mathcal}

\begin{document}

\title[Transmission eigenvalue-free regions]{Transmission eigenvalue-free regions near the real axis. II}

\author[G. Vodev]{Georgi Vodev}

\address {Universit\'e de Nantes, Laboratoire de Math\'ematiques Jean Leray, 2 rue de la Houssini\`ere, BP 92208, 44322 Nantes Cedex 03, France}
\email{Georgi.Vodev@univ-nantes.fr}

\date{}

\begin{abstract} In this paper we extend the results in \cite{kn:V2} to more general domains. More precisely, 
we obtain transmission eigenvalue-free regions for the  
interior transmission problem with one complex-valued refraction index, that is, with a damping term
which does not vanish on the boundary. In particular, we remove the non-trapping condition as well as the strict concavity 
condition from \cite{kn:V2}. Instead, we impose new, more general conditions in terms of the high-frequency behavior of certain cut-off
resolvents associated to exterior problems. 

\end{abstract} 

\maketitle

\setcounter{section}{0}
\section{Introduction}

Our goal in this paper is to extend the recent results on the location of the transmission eigenvalues obtained in \cite{kn:V2}
to more general domains. 
Let $\Omega\subset\mathbb{R}^d$, $d\ge 2$, be a bounded, connected domain with a $C^\infty$ smooth boundary $\Gamma=\partial\Omega$.
In \cite{kn:V2} the following interior transmission problem has been considered:

\begin{equation}\label{eq:1.1}
\left\{
\begin{array}{l}
(\nabla c_1(x)\nabla+\lambda^2n_1(x)+i\lambda m(x))u_1=0\quad \mbox{in}\quad\Omega,\\
(\nabla c_2(x)\nabla+\lambda^2n_2(x))u_2=0\quad \mbox{in}\quad\Omega,\\
u_1=u_2,\,c_1\partial_\nu u_1=c_2\partial_\nu u_2 \quad\mbox{on}\quad\Gamma,
\end{array}
\right.
\end{equation}
where $\lambda\in \mathbb{C}$, $\nu$ denotes the Euclidean unit inner normal to $\Gamma$ and $c_j,n_j,m\in C^\infty(\overline\Omega)$, 
$j=1,2$, are real-valued functions satisfying $c_j(x)>0$, $n_j(x)>0$. We also suppose that either $m(x)\ge 0$ for all $x\in \overline\Omega$
or $m(x)\le 0$ for all $x\in \overline\Omega$. In other words, the function $m$ does not change the sign. 
If the equation (\ref{eq:1.1})
has a non-trivial solution $(u_1,u_2)$ the complex number $\lambda$ is said to be an interior transmission eigenvalue. 
When the function $m$ is not identically zero, it is well-known that there are no real transmission eigenvalues. 
Moreover, Robbiano showed (see Theorem 8 of \cite{kn:R1}) that there are no 
transmission eigenvalues in the region
$$|{\rm Im}\,\lambda|\le C_1e^{-C_2|\lambda|},\quad C_1,C_2>0.$$
When $m\equiv 0$ it is known (e.g. see \cite{kn:NN1}, \cite{kn:NN2}, \cite{kn:LV1}) that the transmission eigenvalues 
form a discreet set in $\mathbb{C}$
with no finite accumulation points under the conditions 
\begin{equation}\label{eq:1.2}
c_1(x)\neq c_2(x) \quad\mbox{on}\quad\Gamma
\end{equation}
 and
\begin{equation}\label{eq:1.3}
c_1(x)n_1(x)\neq c_2(x)n_2(x) \quad\mbox{on}\quad\Gamma.
\end{equation}
This is also true in the case $c_1\equiv c_2\equiv 1$ under the condition (\ref{eq:1.3})
(e.g. see \cite{kn:NN1}, \cite{kn:FN}, \cite{kn:LV2}, \cite{kn:Sy}). 
Note also that under these conditions Weyl asymptotics for the counting function of the transmission eigenvalues 
are obtained in \cite{kn:FN}, \cite{kn:NN2}, \cite{kn:PV}, \cite{kn:R2}. 
We refer to the survey article \cite{kn:CCH} for more results and references in this case.

In the present paper we will be interested in the case when the function $m$ satisfies the condition
\begin{equation}\label{eq:1.4}
m(x)\neq 0\quad\mbox{on}\quad\Gamma.
\end{equation}
The location of the transmission eigenvalues near the real axis has been studied in \cite{kn:V2} under 
(\ref{eq:1.2}) and (\ref{eq:1.4}) as well as the condition that all geodesics associated to the Hamiltonians
$g_j=\frac{c_j(x)}{n_j(x)}|\xi|^2$, $j=1,2$, reach the boundary $\Gamma$ in a finite time. It was also supposed that
$\Gamma$ is $g_2-$ strictly concave (e.g. see \cite{kn:V2} for the precise definition). It is shown in 
\cite{kn:V2} that under these conditions the region
 \begin{equation}\label{eq:1.5}
C_N(|\lambda|+1)^{-N}\le|{\rm Im}\,\lambda|\le C(|\lambda|+1)^{-\ell}
\end{equation}
is free of transmission eigenvalues, $N>1$ being arbitrary, where $\ell\in\{0,1\}$ is defined as follows.
We put $\ell=0$ if the condition 
 \begin{equation}\label{eq:1.6}
\frac{n_1(x)}{c_1(x)}>\frac{n_2(x)}{c_2(x)},\quad c_1(x)<c_2(x),\quad\mbox{on}\quad\Gamma
\end{equation}
is satisfied, and $\ell=1$ otherwise. 

In the present paper we will extend this result to more general situations. To describe them
we first take an open, bounded domain $\mathcal{O}\subset\mathbb{R}^d$ (which may be empty) with a $C^\infty$ smooth boundary $\Gamma_0=\partial\mathcal{O}$,
such that $\mathbb{R}^d\setminus \mathcal{O}$ is connected. Let $\widetilde c_j,\widetilde n_j\in C^\infty(\mathbb{R}^d\setminus \mathcal{O})$, 
$j=1,2$, be real-valued functions satisfying $\widetilde c_j(x)\ge b_0$, $\widetilde n_j(x)\ge b_0$ for some constant $b_0>0$. 
We also suppose that there are 
a constant $R_0\gg 1$ and constants $c_j^0>0$, $n_j^0>0$ such that $\widetilde c_j(x)=c_j^0$, $\widetilde n_j(x)=n_j^0$ for $|x|\ge R_0+1$. 
Denote by $\widetilde G_j$
the selfadjoint Dirichlet realization of the operator $-\widetilde n_j(x)^{-1}\nabla\widetilde c_j(x)\nabla$ on the Hilbert space 
$\widetilde H_j=L^2(\mathbb{R}^d\setminus \mathcal{O};\widetilde n_j(x)dx)$. We define the outgoing and incoming resolvents, $\mathcal{R}_j^+(\lambda)$
and $\mathcal{R}_j^-(\lambda)$, so that 
$$\mathcal{R}_j^\pm(\lambda)=(\widetilde G_j-\lambda^2)^{-1}:\widetilde H_j\to \widetilde H_j$$
is bounded for $\pm{\rm Im}\,\lambda>0$. Clearly, we have $\mathcal{R}_j^+(\lambda)^*=\mathcal{R}_j^-(\overline\lambda)$.
Let $\chi\in C^\infty(\mathbb{R}^d)$
be a real-valued function of compact support such that $\chi(x)=1$ for $|x|\le R_0+2$. 
It is well-known that the cutoff resolvents $\chi\mathcal{R}_j^+(\lambda)\chi$
and $\chi\mathcal{R}_j^-(\lambda)\chi$ extend through the real axis to meromorphic operator-valued functions with no poles on the real axis. 
Moreover, we have 
$$\left(\chi\mathcal{R}_j^+(\lambda)\chi\right)^*=\chi\mathcal{R}_j^-(\lambda)\chi$$
 for real $\lambda$.
We suppose that the cut-off resolvent of $\widetilde G_j$ satisfies the following
high-frequency bound on the real axis:
\begin{equation}\label{eq:1.7}
\left\|\partial_x^\alpha\chi\mathcal{R}_j^\pm(\lambda)\chi\partial_x^\beta\right\|_{\widetilde H_j\to \widetilde H_j}\lesssim \lambda^{|\alpha|+|\beta|-1}\tau_j(\lambda),\quad \lambda\gg 1,
\end{equation}
for all multi-indices $\alpha$ and $\beta$ such that $|\alpha|+|\beta|\le 2$, where the function $\tau_j$ is non-decreasing
and satisfies
\begin{equation}\label{eq:1.8}
1\le \tau_j(\lambda)\le \lambda^q,
\end{equation}
with some constant $q\ge 0$. 
Let $\Omega_\sharp\subset\mathbb{R}^d$ be a bounded, connected domain with a $C^\infty$ smooth boundary $\Gamma=\partial\Omega_\sharp$, such that 
$$\overline{\mathcal{O}}\subset \Omega_\sharp\subset \{x\in \mathbb{R}^d: |x|\le R_0\}.$$
Then $\Omega=\Omega_\sharp\setminus \overline{\mathcal{O}}$ is a bounded, connected domain with a $C^\infty$ smooth boundary 
$\partial\Omega=\Gamma\cup \Gamma_0$ and $\Gamma\cap \Gamma_0=\emptyset$. We denote by $c_j$, $n_j$ the restrictions of the functions
$\widetilde c_j$, $\widetilde n_j$ on $\Omega$. In the present paper we will consider the following interior transmission problem:

\begin{equation}\label{eq:1.9}
\left\{
\begin{array}{l}
(\nabla c_1(x)\nabla+\lambda^2n_1(x)+i\lambda m(x))u_1=0\quad \mbox{in}\quad\Omega,\\
(\nabla c_2(x)\nabla+\lambda^2n_2(x))u_2=0\quad \mbox{in}\quad\Omega,\\
u_1=u_2,\,c_1\partial_\nu u_1=c_2\partial_\nu u_2 \quad\mbox{on}\quad\Gamma,\\
u_1=u_2=0\quad\mbox{on}\quad\Gamma_0,
\end{array}
\right.
\end{equation}
where the function $m$ is as above. We will be interested in obtaining eigenvalue-free regions near the 
real axis under the above conditions. 
Our first result is the following

\begin{Theorem} \label{1.1}
Suppose that the conditions (\ref{eq:1.4}) and (\ref{eq:1.7}) are satisfied with $j=1$.  
 Then there exists a constant $C>0$ such that there are no transmission eigenvalues in the region
 \begin{equation}\label{eq:1.10}
-C\tau_1(|\lambda|)^{-2}\le {\rm sign}(m)\,{\rm Im}\,\lambda\le 0,\quad |\lambda|>1.
\end{equation}
 \end{Theorem}

This theorem is relatively easy to prove. It follows from the a priori estimates for the solution of the first equation in
(\ref{eq:1.9}) obtained in Section 3 together with the Green formula applied to the second equation in
(\ref{eq:1.9}). Therefore we do not need to impose any conditions on the coefficients $c_2$ and $n_2$.

In the half-plane $\{{\rm sign}(m)\,{\rm Im}\,\lambda>0\}$, however, the situation is much more complicated. 
Therefore, our main goal in the
present paper is to obtain eigenvalue-free regions also in this  
half-plane under the conditions
(\ref{eq:1.2}), (\ref{eq:1.4}), (\ref{eq:1.7}) and (\ref{eq:1.8}). Define $\ell_1,\ell_2\in\{0,1\}$ as follows. We put 
$\ell_1=0$ if the condition (\ref{eq:1.6}) is satisfied and the function $\tau_1$ satisfies
\begin{equation}\label{eq:1.11}
\lambda^{-1/4}\tau_1(\lambda)\to 0\quad\mbox{as}\quad \lambda\to \infty,
\end{equation}
and we put $\ell_1=1$ otherwise. We put 
$\ell_2=0$ if $\Gamma$ is $g_2-$ strictly concave and the function $\tau_2$ satisfies
\begin{equation}\label{eq:1.12}
\lambda^{-1/2}\tau_2(\lambda)\to 0\quad\mbox{as}\quad \lambda\to \infty,
\end{equation}
and we put $\ell_2=1$ otherwise. Given a parameter $C>0$ we denote by $\Lambda(C)\subset\mathbb{C}$ the region 
$$|{\rm Im}\,\lambda|\le C|\lambda|^{-\ell_1-\ell_2}\tau_1(|\lambda|)^{-10+4\ell_1}\tau_2(|\lambda|)^{-2},\quad |\lambda|>1.$$
Denote by $G_2$
the selfadjoint Dirichlet realization of the operator $-n_2(x)^{-1}\nabla c_2(x)\nabla$ on the Hilbert space 
$H_2=L^2(\Omega;n_2(x)dx)$. Clearly, the operator $G_2$ is non-negative. It is also well-known that 
the spectrum of $G_2$ is discreet consisting of infinitely many non-negative eigenvalues. 
Given $\lambda\in\mathbb{C}^\pm:=\{\lambda\in\mathbb{C}:\pm{\rm Re}\,\lambda>0\}$, set 
$$\theta_\pm(\lambda)={\rm dist}(\lambda,{\rm spec}\,\pm\sqrt{G_2})=\min_{\mu_k\in {\rm spec}\,\sqrt{G_2}}|\lambda\mp\mu_k|
\ge |{\rm Im}\,\lambda|.$$ 
Given parameters $C,N>0$ we denote by $\mathcal{L}_N^\pm(C)\subset\mathbb{C}^\pm$ the region 
$$\theta_\pm(\lambda)\le C|\lambda|^{-N},\quad |\lambda|>1.$$
Let $\mathcal{L}_N(C)=\mathcal{L}_N^+(C)\cup\mathcal{L}_N^-(C)$. 
We are now ready to state our main result.

\begin{Theorem} \label{1.2}
Suppose that the conditions (\ref{eq:1.2}), (\ref{eq:1.4}), (\ref{eq:1.7}) and (\ref{eq:1.8}) are satisfied.  
 Then for every $N>12q+2$ there exist constants $C, C_N>0$,
$C$ being independent of $N$, such that there are no transmission eigenvalues in the region
$\Lambda(C)\setminus \mathcal{L}_N(C_N)$. 
\end{Theorem}

\begin{rem} \label{1.3} 
Obviously, if $\lambda$ is a transmission eigenvalue associated to the problem (\ref{eq:1.9}), then 
$-\lambda$ is a transmission eigenvalue associated to the problem (\ref{eq:1.9}) with the function $m$ replaced by $-m$.
Therefore, it suffices to prove the above theorems in $\mathbb{C}^+$, only.
 \end{rem}

\begin{rem} \label{1.4} It follows from Theorem \ref{1.2} that the transmission eigenvalues in the region
$\Lambda(C)\cap\mathbb{C}^\pm$ are either finitely many or they are infinitely many and 
asymptotically very close to the eigenvalues of the
self-adjoint operator $\pm\sqrt{G_2}$. 
\end{rem}

\begin{rem} \label{1.5} It is well-known (e.g. see \cite{kn:B1}) that the cutoff resolvents always satisfy the bound (\ref{eq:1.7})
with $\tau_j=e^{c\lambda}$, $c>0$ being some constant, and that this bound cannot be improved without imposing extra conditions. 
Therefore, imposing the polynomial bound (\ref{eq:1.8}) is important in the proof of Theorem \ref{1.2}. 
Indeed, our arguments do not work anymore without this condition. In particular, it guarantees that $\mathcal{L}_N(C_N)\subset\Lambda(C)$,
provided $N>12q+2$. 
The condition (\ref{eq:1.8}), however, is not necessary in the proof of Theorem \ref{1.1}. 
\end{rem}

\begin{rem} \label{1.6} When $\mathcal{O}=\emptyset$ and the coefficients $\widetilde c_j, \widetilde n_j$, $j=1,2$, are constants,
then $\mathcal{R}_j^+(\lambda)$ and $\mathcal{R}_j^-(\lambda)$ are the outgoing and incoming resolvents of the operator
$-\frac{\widetilde c_j}{\widetilde n_j}\Delta$, where $\Delta$ is the Euclidean Laplacian on $\mathbb{R}^d$.
In this case it is well-known that the bound (\ref{eq:1.7}) holds 
with $\tau_j=1$.  Such a bound still holds for non-constant coefficients but one has to assume a non-trapping condition,
namely that all geodesics associated to the Hamiltonians $\frac{\widetilde c_j(x)}{\widetilde n_j(x)}|\xi|^2$ escape to the infinity.
This conclusion remains valid when $\mathcal{O}\neq\emptyset$ as well, but in this case one has to impose the non-trapping condition
on the broken geodesics. The simplest example of such a non-trapping obstacle $\mathcal{O}$ is the strictly convex one, provided 
the coefficients $\widetilde c_j, \widetilde n_j$, $j=1,2$, are supposed constants. 
\end{rem}

\begin{rem} \label{1.7} Our result also applies to trapping obstacles $\mathcal{O}$ which may have periodic broken geodesics,
provided these geodesics are of hyperbolic type. A typical example is an obstacle consisting of several strictly convex bodies and
the coefficients are supposed constants (in which case the geodesics are just rays). More precisely, let 
$\mathcal{O}=\cup_{j=1}^J\mathcal{O}_j$, $J\ge 2$, where each $\mathcal{O}_j$ is strictly convex and 
$\overline{\mathcal{O}_i}\cap\overline{\mathcal{O}_j}=\emptyset$ if $i\neq j$. When $J=2$ there is only one periodic
ray between $\mathcal{O}_1$ and $\mathcal{O}_2$ (which is of hyperbolic type), 
while when $J\ge 3$ there are infinitely many. Therefore, in this case
one needs to impose some conditions in order to be able to get a nice resolvent bound on the real axis. 
The first one is the Ikawa no-eclipse condition: if $i$, $j$ and $k$ are all different, then
$\overline{\mathcal{O}_k}$ does not intersect the convex hull of $\overline{\mathcal{O}_i}\cup\overline{\mathcal{O}_j}$. 
The second one is a dynamical assumption involving the topological presure of the billiard flow (see \cite{kn:I} for
the precise definition). Under these conditions Ikawa \cite{kn:I} obtained that the cut-off resolvent extends analytically to
a strip $|{\rm Im}\,\lambda|\le C$ and its norm is polynomially bounded there. Using this Burq \cite{kn:B2} applied the Fragm\`en-Lindel\"of
principle to conclude that on the real axis the norm of the cut-off resolvent is logarithmically bounded. Therefore,
under the Ikawa conditions mentioned above, in this case we have the bound (\ref{eq:1.7}) fulfilled with $\tau_j=\log\lambda$. 
\end{rem}

The eigenvaue-free region in Theorem \ref{1.2} is more precise than the region (\ref{eq:1.5}) obtained in \cite{kn:V2}
since $\mathcal{L}_N(C_N)\subset\{|{\rm Im}\,\lambda|\le C_N|\lambda|^{-N}\}$. 
 This is due to the observation that the parametrix of the interior Dirichlet-to-Neumann map (without damping) 
in the elliptic region is valid in $\mathbb{C}\setminus\mathcal{L}_N(C_N)$ (see Section 6). Note that this parametrix plays a
crucial role in our proof. Therefore, studying the transmission eigenvalues in $\mathcal{L}_N(C_N)$, that is,
close to the eigenvalues of $\pm\sqrt{G_2}$, requires a different approach.
Most probably the region $\mathcal{L}_N(C_N)$ is also free of transmission eigenvalues but proving this remains an open problem.

To prove Theorem \ref{1.2} we follow the same strategy as in \cite{kn:V2} with some important modifications. The first one is, 
as mentioned above, that we extend the domain of validity of the parametrix of the interior Dirichlet-to-Neumann map 
in the elliptic region (see Theorem \ref{6.1}).
Secondly, we obtain more general a priori estimates for the solutions of the interior Helmholtz equation without
damping (see Theorem \ref{2.1}). In particular, the boundary $\Gamma$ is not necessairily supposed strictly concave. This is done by using
the so-called jump formula that allows to express the solution by the term in the right-hand side and the restriction on $\Gamma$
of the solution and its normal derivative (see Lemma \ref{2.5}). The a priori estimates for the solutions of the interior 
Helmholtz equation with damping are also more general than those obtained in \cite{kn:V2} (see Theorem \ref{3.1}).

The proof of Theorems \ref{1.1} and \ref{1.2} is carried out in Section 7. Let us sketch our approach
to the proof of Theorem \ref{1.2} which is much more complicated. 
We need to show that if $\lambda$ belongs to the eigenvalue-free regions, 
then the solution $(u_1,u_2)$ to the equation 
(\ref{eq:1.9}) is identically zero. In fact, it suffices to show that the function $f=u_1|_\Gamma=u_2|_\Gamma$ is identically zero, 
which would imply that $u_1$ and $u_2$ are identically zero. Clearly, if $\lambda$ is a transmission eigenvalue, we have
 $T(\lambda)f\equiv 0$ with $f$ not identically zero, where 
$$T(\lambda)=c_1{\cal N}_1(\lambda,m)-c_2{\cal N}_2(\lambda),$$ 
${\cal N}_1(\lambda,m)$ and ${\cal N}_2(\lambda)$ being the corresponding interior Dirichlet-to-Neumann maps (see Section 5).
In other words, we have to show that if $\lambda$ belongs to the eigenvalue-free region of Theorem \ref{1.2} and
$T(\lambda)f\equiv 0$, then $f\equiv 0$. To this end, we use the parametrix of the interior Dirichlet-to-Neumann maps in the
corresponding elliptic regions in the form of $h-\Psi$DOs (see Section 6). 
This allows us to build a parametrix for the operator $T(\lambda)$
in the intersection of the two elliptic regions and to compute its principal symbol. The assumption (\ref{eq:1.2}) is crucial
since it implies that this principal symbol is an elliptic one in the deep elliptic region belonging to the class $S^1(\Gamma)$ 
(see Section 4 for the definition). This fact in turm implies that the parametrix of $T(\lambda)$ sends the Sobolev space
$H^1(\Gamma)$ into $L^2(\Gamma)$ and it is invertible in the deep elliptic region. That is what we use in the proof of Lemma
\ref{7.2}. The assumption (\ref{eq:1.6}) allows to extend the invertibility of the parametrix to the whole elliptic region associated to
the first equation in (\ref{eq:1.9}). This is used in the proof of Lemma \ref{7.3}. Finally, we combine Lemmas \ref{7.2} and \ref{7.3}
with the a priori estimates from the previous sections to obtain the desired eiganvalue-free region.

\section{Study of the interior Helmholtz equation without damping}

Throughout this paper $\|\cdot\|$, $\|\cdot\|_1$, $\|\cdot\|_0$ and $\|\cdot\|_{1,0}$ will denote the norms in $L^2(\Omega)$,
$H^1(\Omega)$, $L^2(\Gamma)$ and $H^1(\Gamma)$, respectively, where the norms in $H^1(\Omega)$ and $H^1(\Gamma)$ are the semiclassical ones, that is, 
$$\|u\|_1^2:=\sum_{0\le|\alpha|\le 1}\left\|(h\partial_x)^\alpha u\right\|^2,$$
$$\|u\|_{1,0}^2:=\sum_{0\le|\alpha|\le 1}\left\|(h\partial_x)^\alpha u\right\|_0^2,$$
where $0<h\ll 1$ is a semiclassical parameter to be fixed below. 
Also, $\langle\cdot,\cdot\rangle$ and $\langle\cdot,\cdot\rangle_0$ will denote the scalar products in 
$L^2(\Omega)$ and $L^2(\Gamma)$, respectively. 

In this section we consider the equation
\begin{equation}\label{eq:2.1}
\left\{
\begin{array}{l}
(\nabla c(x)\nabla+\lambda^2n(x))u=\lambda v\quad \mbox{in}\quad\Omega,\\
u=f\quad\mbox{on}\quad\Gamma,\\
u=0\quad\mbox{on}\quad\Gamma_0,
\end{array}
\right.
\end{equation}
where $\lambda\in\mathbb{C}$, ${\rm Re}\,\lambda\gg 1$, and $c,n\in C^\infty(\overline\Omega)$ are real-valued functions satisfying $c(x)>0$, $n(x)>0$ for all $x\in \overline\Omega$. 
Let $\widetilde c,\widetilde n\in C^\infty(\mathbb{R}^d\setminus \mathcal{O})$ be real-valued functions satisfying 
$\widetilde c(x)\ge b_0$, $\widetilde n(x)\ge b_0$, $b_0>0$, such that $\widetilde c=c$, $\widetilde n=n$ in $\Omega$. 
We also suppose that there are 
a constant $R_0\gg 1$ and constants $c^0>0$, $n^0>0$ such that $\overline\Omega\subset \{x\in \mathbb{R}^d: |x|\le R_0\}$ 
and $\widetilde c(x)=c^0$, $\widetilde n(x)=n^0$ for $|x|\ge R_0+1$. 
Denote by $\widetilde G$
the selfadjoint Dirichlet realization of the operator $-\widetilde n(x)^{-1}\nabla\widetilde c(x)\nabla$ on the Hilbert space 
$\widetilde H=L^2(\mathbb{R}^d\setminus \mathcal{O};\widetilde n(x)dx)$. 
We define the outgoing and incoming resolvents, $\mathcal{R}^+(\lambda)$
and $\mathcal{R}^-(\lambda)$, so that 
$$\mathcal{R}^\pm(\lambda)=(\widetilde G-\lambda^2)^{-1}:\widetilde H\to \widetilde H$$
is bounded for $\pm{\rm Im}\,\lambda>0$. Clearly, we have $\mathcal{R}^+(\lambda)^*=\mathcal{R}^-(\overline\lambda)$.
Let $\chi\in C^\infty(\mathbb{R}^d)$
be a real-valued function of compact support such that $\chi(x)=1$ for $|x|\le R_0+2$. 
As mentioned in the introduction, the cutoff resolvents $\chi\mathcal{R}^+(\lambda)\chi$
and $\chi\mathcal{R}^-(\lambda)\chi$ extend through the real axis to meromorphic operator-valued functions with no poles on the real axis,  
and we have $\left(\chi\mathcal{R}^+(\lambda)\chi\right)^*=\chi\mathcal{R}^-(\lambda)\chi$ for real $\lambda$.
We suppose that the cut-off resolvent of $\widetilde G$ satisfies the bound 
\begin{equation}\label{eq:2.2}
\left\|\partial_x^\alpha\chi\mathcal{R}^\pm(\lambda)\chi\partial_x^\beta\right\|_{\widetilde H\to \widetilde H}\lesssim \lambda^{|\alpha|+|\beta|-1}\tau(\lambda),\quad \lambda\gg 1,
\end{equation}
for all multi-indices $\alpha$ and $\beta$ such that $|\alpha|+|\beta|\le 2$, where the function $\tau\ge 1$ is non-decreasing,
that is, $\tau(\lambda_1)\le \tau(\lambda_2)$ if $\lambda_1\le\lambda_2$.

In $\Omega$ we define the Hamiltonian $g=\frac{c(x)}{n(x)}|\xi|^2$. 
We now define $\ell\in\{0,1\}$ as follows. 
We put $\ell=0$ if $\Gamma$ is $g-$ strictly concave and $\tau$ satisfies
\begin{equation}\label{eq:2.3}
\lambda^{-1/2}\tau(\lambda)\to 0\quad\mbox{as}\quad \lambda\to \infty,
\end{equation}
and we put $\ell=1$ otherwise. We introduce the semiclassical parameter $h=({\rm Re}\,\lambda)^{-1}$. 
Given a parameter $0<\epsilon\ll 1$, independent of $\lambda$, set 
$\Omega_\epsilon=\{x\in\Omega:{\rm dist}(x,\Gamma)<\epsilon\}$.  
One of our goals in this section is to prove the following

\begin{Theorem} \label{2.1} 
Suppose that the condition (\ref{eq:2.2}) is fulfilled. Let $u\in H^2(\Omega)$ satisfy equation (\ref{eq:2.1}) and set 
 $\omega=h\partial_\nu u|_\Gamma$. Then there are constants $C,\lambda_0>0$ such that for all $\lambda\in\mathbb{C}$ such that 
 \begin{equation}\label{eq:2.4}
 |{\rm Im}\,\lambda|\le C\tau(|\lambda|)^{-1},\quad {\rm Re}\,\lambda\ge\lambda_0, 
 \end{equation}
 we have the estimates
\begin{equation}\label{eq:2.5}
\|u\|\lesssim \tau(|\lambda|)\left(\|v\|+\|u\|_{L^2(\Omega_\epsilon)}\right),
\end{equation}
\begin{equation}\label{eq:2.6}
\|u\|_1\lesssim \tau(|\lambda|)\left(\|v\|+|\lambda|^{\ell/2}\|f\|_0+|\lambda|^{\ell/2}\|\omega\|_0\right).
\end{equation}
\end{Theorem}

{\it Proof.} It is easy to see that it suffices to prove the estimates (\ref{eq:2.5}) and (\ref{eq:2.6}) for real $\lambda\ge\lambda_0$.
Indeed, if $u$ satisfies equation (\ref{eq:2.1}) with complex $\lambda$, then $u$ satisfies equation (\ref{eq:2.1}) with
$\lambda$ replaced by ${\rm Re\,\lambda}$ and $v$ replaced by 
$$\widetilde v=\frac{\lambda}{{\rm Re\,\lambda}}v+\frac{{\rm Im\,\lambda}}{{\rm Re\,\lambda}}({\rm Im\,\lambda}-2i{\rm Re\,\lambda})nu.$$
Therefore, by (\ref{eq:2.6}) applied with $\lambda$ replaced by ${\rm Re\,\lambda}$ and $v$ replaced by $\widetilde v$ we get
\begin{equation}\label{eq:2.7}
\|u\|_1\lesssim \tau(|\lambda|)\left(\|v\|+|\lambda|^{\ell/2}\|f\|_0+|\lambda|^{\ell/2}\|\omega\|_0\right)
+|{\rm Im}\,\lambda|\tau(|\lambda|)\|u\|.
\end{equation}
If $\lambda$ satisfies (\ref{eq:2.4}), taking $C$ properly we can absorb the last term in the right-hand side of (\ref{eq:2.7}). 
Therefore we can conclude that
(\ref{eq:2.6}) also holds for $\lambda$ satisfying (\ref{eq:2.4}). Clearly, a similar analysis applies to the estimate (\ref{eq:2.5}).
Thus, in what follows we will prove the estimates (\ref{eq:2.5}) and (\ref{eq:2.6}) 
for real $\lambda\ge\lambda_0$. Then we have $h=\lambda^{-1}$.

Let $\phi,\phi_1\in C_0^\infty(\mathbb{R}^d)$ be independent of $\lambda$ and such that $\phi=1$ in 
$(\Omega\cup\mathcal{O})\setminus\Omega_{\epsilon/3}$, supp$\,\phi\subset(\Omega\cup\mathcal{O})\setminus\Omega_{\epsilon/4}$,
and $\phi_1=1$ in
$(\Omega\cup\mathcal{O})\setminus\Omega_{\epsilon}$, supp$\,\phi_1\subset(\Omega\cup\mathcal{O})\setminus\Omega_{\epsilon/2}$. 

\begin{lemma} \label{2.2} 
The solution $u$ of the equation (\ref{eq:2.1}) satisfies the formula
\begin{equation}\label{eq:2.8}
\phi u=-\chi\mathcal{R}^\pm(\lambda)\chi n^{-1}\left(\lambda \phi v+[\nabla c\nabla,\phi](1-\phi_1)u\right)
\end{equation}
for real $\lambda$.
\end{lemma}

{\it Proof.} Clearly, the function $\phi u$ satisfies the equation

\begin{equation}\label{eq:2.9}
\left\{
\begin{array}{l}
(\nabla c(x)\nabla+\lambda^2n(x))\phi u=U\quad \mbox{in}\quad\mathbb{R}^d\setminus \mathcal{O},\\
\phi u=0\quad\mbox{on}\quad\Gamma_0,
\end{array}
\right.
\end{equation}
where $U\in \widetilde H$ is given by
$$U=\lambda \phi v+[\nabla c(x)\nabla,\phi]u=\lambda \phi v+[\nabla c(x)\nabla,\phi](1-\phi_1)u.$$
Given a parameter $0<k\ll 1$, set
$$U_k^\pm=(\nabla c(x)\nabla+(\lambda\pm ik)^2n(x))\phi u.$$
We have 
$$\left\|U-U_k^\pm\right\|_{\widetilde H}\lesssim k|\lambda|\|\phi u\|_{\widetilde H},$$
which implies
$$\lim_{k\to 0}\left\|U-U_k^\pm\right\|_{\widetilde H}=0.$$
On the other hand, since $\chi=1$ on supp$\,\phi$, we have
$$\phi u=-\chi\mathcal{R}^\pm(\lambda\pm ik)\chi n^{-1}U_k^\pm.$$
Thus, taking the limit $k\to 0$ we arrive at the formula
$$\phi u=-\chi\mathcal{R}^\pm(\lambda)\chi n^{-1}U,$$
as desired.
\eproof

By (\ref{eq:2.2}) with $\alpha=0$, $|\beta|=1$, and (\ref{eq:2.8}), we obtain
\begin{equation}\label{eq:2.10}
\|\phi u\|\lesssim \tau(\lambda)\|v\|+\tau(\lambda)\|(1-\phi_1)u\|,
\end{equation}
which clearly implies (\ref{eq:2.5}). 
When $\ell=0$ the estimate (\ref{eq:2.6}) for real $\lambda\gg 1$ follows from  (\ref{eq:2.10}) and the following

\begin{prop} \label{2.3} 
Let $\Gamma$ be $g-$strictly concave. Then, for a suitable choice of the function
$\phi_1$ and the parameter $\epsilon$, both independent of $\lambda$, we have the estimate
\begin{equation}\label{eq:2.11}
\|(1-\phi_1)u\|\lesssim \|v\|+\|f\|_0+\|\omega\|_0+h^{1/2}\|u\|_1.
\end{equation}
\end{prop}

Indeed, combining the estimates (\ref{eq:2.10}) and (\ref{eq:2.11}) leads to
\begin{equation}\label{eq:2.12}
\|u\|\lesssim \tau(\lambda)\left(\|v\|+\|f\|_0+\|\omega\|_0\right)+\lambda^{-1/2}\tau(\lambda)\|u\|_1. 
\end{equation}
On the other hand, by the Green formula we have
\begin{equation}\label{eq:2.13}
\langle\lambda^2nu-\lambda v,u\rangle=
\langle -\nabla c\nabla u,u\rangle=\int_\Omega c|\nabla u|^2+h^{-1}\langle c\omega,f\rangle_0,
\end{equation}
which clearly still holds for complex $\lambda$. 
Taking the real part of (\ref{eq:2.13}) leads to the estimate
\begin{equation}\label{eq:2.14}
\|u\|_1\lesssim h\|v\|+\|u\|+h^{1/2}\|f\|_0^{1/2}\|\omega\|_0^{1/2}.
\end{equation}
Hence
\begin{equation}\label{eq:2.15}
\|u\|_1\lesssim \|u\|+\|v\|+\|f\|_0+\|\omega\|_0.
\end{equation}
We now combine 
the estimates (\ref{eq:2.12}) and (\ref{eq:2.15}) and use (\ref{eq:2.3}) to absorb the term $\lambda^{-1/2}\tau(\lambda)\|u\|_1$
by taking $\lambda$ big enough. Clearly, this leads to the estimate (\ref{eq:2.6}) in this case.
Note that the above proposition is in fact Proposition 2.2 of \cite{kn:CPV} and we refer the reader 
to Section 2 of \cite{kn:CPV} for the proof.

Consider now the case $\ell=1$. Let $\gamma:H^s(\Omega)\to H^{s-1/2}(\Gamma)$, $s>1/2$, denote the restriction on $\Gamma$ and let 
 $\gamma^*:H^{-s+1/2}(\Gamma)\to H^{-s}(\Omega)$ be its adjoint. 
We will need the following

\begin{lemma} \label{2.4} 
The restriction on the boundary satisfies the estimate
\begin{equation}\label{eq:2.16}
\|\gamma u\|_0\lesssim \lambda^{1/2}\left(\varepsilon\|u\|_1+\varepsilon^{-1}\|u\|\right)
\end{equation}
for every $0<\varepsilon\le 1$. 
\end{lemma}

{\it Proof.} Clearly, it suffices to prove (\ref{eq:2.16}) locally near $\Gamma$. 
Let ${\cal V}\subset\mathbb{R}^d$ be a small open domain such that ${\cal V}^0:={\cal V}\cap\Gamma\neq\emptyset$. 
Let $(x_1,x')\in {\cal V}^+:={\cal V}\cap\Omega$, $0<x_1\ll 1$, $x'=(x_2,...,x_d)\in{\cal V}^0$, 
be the local normal geodesic coordinates near the boundary. 
Let also ${\cal V}_1\subset{\cal V}$ be a small open domain such that ${\cal V}_1^0:={\cal V_1}\cap\Gamma\neq\emptyset$. 
Choose a function $\psi\in C_0^\infty({\cal V})$, $0\le\psi\le 1$, such that $\psi=1$ on ${\cal V}_1$. 
Set $u^\flat:=\psi(1-\phi)u$, the function $\phi$ being as above, and $\mathcal{D}_{x_1}=-ih\partial_{x_1}$. We have
 $$-\frac{d}{dx_1}\|u^\flat(x_1,\cdot)\|_0^2=-2{\rm Re}\langle u^\flat(x_1,\cdot),\partial_{x_1}u^\flat(x_1,\cdot)\rangle_0$$ 
 $$\le h^{-1}\varepsilon^{-2}\|u^\flat(x_1,\cdot)\|_0^2+h^{-1}\varepsilon^2\|\mathcal{D}_{x_1}u^\flat(x_1,\cdot)\|_0^2.$$
 Hence
 $$\|\gamma\psi\gamma u\|_0^2=\|u^\flat(0,\cdot)\|_0^2=-\int_0^\epsilon \frac{d}{dx_1}\|u^\flat(x_1,\cdot)\|_0^2dx_1$$
 $$\le h^{-1}\varepsilon^{-2}\|u^\flat\|^2+h^{-1}\varepsilon^2\|\mathcal{D}_{x_1}u^\flat\|^2\lesssim h^{-1}\varepsilon^{-2}\|u\|^2+h^{-1}\varepsilon^2\|u\|_1^2.$$
  Since $\Gamma$ is compact, there exist a finite number of smooth functions $\psi_i$, $0\le\psi_i\le 1$, $i=1,...,I,$ such that 
$1=\sum_{i=1}^I\psi_i$ and the above estimate holds with $\gamma\psi$ replaced by each $\psi_i$. Therefore, 
 the estimate (\ref{eq:2.16}) is obtained by summing up all such estimates.
\eproof

To prove (\ref{eq:2.6}) in this case we will express the solution $u$
 in terms of the resolvent $\mathcal{R}^\pm(\lambda)$ and the functions $v$, $f$ and $\omega$ by using the so-called jump formula. 
 It can be derived from (\ref{eq:2.8}) by letting $\epsilon\to 0$. We have the following
 
 \begin{lemma} \label{2.5} 
 There exists a first-order differential operator $Q$ near $\Gamma$ so that 
 the solution to the equation (\ref{eq:2.1}) satisfies the formula
 \begin{equation}\label{eq:2.17}
\phi_0 u=-\lambda \chi\mathcal{R}^\pm(\lambda)\chi n^{-1}\left(\phi_0 v+hQ\gamma^*f+\gamma^*c_0\omega\right),
\end{equation}
for real $\lambda$, where $c_0=\gamma c$ and $\phi_0$ denotes the characteristic function of $\Omega\cup\mathcal{O}$.
 \end{lemma}
 
 {\it Proof.} We take a family of smooth real-valued functions $\phi_\epsilon\to \phi_0$ as 
 $\epsilon\to 0$ for which (\ref{eq:2.8}) holds. 
 We will derive from the Green formula that
 \begin{equation}\label{eq:2.18}
 \langle[\nabla c\nabla,\phi_\epsilon]u,w\rangle\to 
 \langle c_0\gamma\partial_\nu u,\gamma w\rangle_0+\langle \gamma u,\gamma\widetilde Q w\rangle_0
 \end{equation}
  for every $w\in H^1(\Omega)$ such that $w=0$ on $\Gamma_0$, where $\widetilde Q$ is a first-order differential operator. 
 Indeed, by the Green formula we have the identity
 \begin{equation}\label{eq:2.19}
 \langle -\nabla c\nabla u,w\rangle=\langle c\nabla u,\nabla w\rangle,
 \end{equation} 
 provided either $w=0$ or $\partial_\nu u=0$ on $\partial\Omega$.
 Using (\ref{eq:2.19}) we obtain
 \begin{equation}\label{eq:2.20}
 \langle [\nabla c\nabla,\phi_\epsilon] u,w\rangle=\langle c\nabla u,[\nabla,\phi_\epsilon]w\rangle
 -\langle c[\nabla,\phi_\epsilon]u,\nabla w\rangle.
 \end{equation} 
 On the other hand, near $\Gamma$ we have $\nabla=\partial_\nu+\mathcal{Q}$, where $\mathcal{Q}$ is a tangential first-order 
 differential operator. Therefore,
 $\gamma\nabla u=\gamma\partial_\nu u+\widetilde{\mathcal{Q}}\gamma u$, where $\widetilde{\mathcal{Q}}$ is the first-order
 differential operator obtained by restricting the coefficients of $\mathcal{Q}$ on $\Gamma$. We also have in the sense of distributions
 $$[\nabla,\phi_\epsilon]\to [\partial_\nu,\phi_\epsilon]\to [\partial_\nu,\phi_0]=\delta_\Gamma,$$
 where $\delta_\Gamma$ denotes the Dirac delta function on $\Gamma$. Thus, by (\ref{eq:2.20}) we get
 $$ \langle[\nabla c\nabla,\phi_\epsilon]u,w\rangle\to 
 \langle c_0\gamma\partial_\nu u,\gamma w\rangle_0-\langle c_0\gamma u,\gamma\partial_\nu w\rangle_0$$
  $$+\langle c_0\widetilde{\mathcal{Q}}\gamma u,\gamma w\rangle_0-\langle c_0\gamma u,\widetilde{\mathcal{Q}}\gamma w\rangle_0$$
  $$=\langle c_0\gamma\partial_\nu u,\gamma w\rangle_0-\langle c_0\gamma u,\gamma\partial_\nu w\rangle_0$$
  $$+\langle\gamma u,(\widetilde{\mathcal{Q}}^*c_0-c_0\widetilde{\mathcal{Q}})\gamma w\rangle_0,$$
  which implies (\ref{eq:2.18}) with 
  $$\widetilde Q=-c_0\partial_\nu+\widetilde{\mathcal{Q}}^*c_0-c_0\widetilde{\mathcal{Q}}.$$
  Given any $w_1\in L^2(\Omega)$, by (\ref{eq:2.8}) we have
  $$-\langle \phi_\epsilon u,w_1\rangle=\lambda \langle\chi\mathcal{R}^\pm(\lambda)\chi n^{-1}\phi_\epsilon v,w_1\rangle
  +\langle[\nabla c\nabla,\phi_\epsilon]u,n^{-1}\chi\mathcal{R}^\mp(\lambda)\chi w_1\rangle.$$
  Letting $\epsilon\to 0$ and using (\ref{eq:2.18}) with $w=n^{-1}\chi\mathcal{R}^\mp(\lambda)\chi w_1$, we get
  \begin{equation}\label{eq:2.21}
  -\langle \phi_0 u,w_1\rangle=\lambda \langle\chi\mathcal{R}^\pm(\lambda)\chi n^{-1}\phi_0 v,w_1\rangle$$
  $$+\langle c_0\gamma\partial_\nu u,\gamma n^{-1}\chi\mathcal{R}^\mp(\lambda)\chi w_1\rangle_0
  +\langle \gamma u,\gamma\widetilde Q n^{-1}\chi\mathcal{R}^\mp(\lambda)\chi w_1\rangle_0,
  \end{equation} 
 which clearly implies (\ref{eq:2.17}) with $Q=\widetilde Q^*$.
\eproof

We will now use the form (\ref{eq:2.21}) of the formula (\ref{eq:2.17}) together with Lemma \ref{2.4} applied with $\varepsilon=1$
and assumption (\ref{eq:2.2}). 
We get
$$|\langle \phi_0 u,w_1\rangle|\le \lambda \|\chi\mathcal{R}^\pm(\lambda)\chi n^{-1}\phi_0 v\|\|w_1\|$$
$$+\lambda^{3/2} \|\chi\mathcal{R}^\mp(\lambda)\chi w_1\|_1\|\omega\|_0
+\lambda^{3/2}\|hQ\chi\mathcal{R}^\mp(\lambda)\chi w_1\|_1\|f\|_0$$
$$\lesssim\tau(\lambda)\|w_1\|(\|v\|+\lambda^{1/2}\|f\|_0+\lambda^{1/2}\|\omega\|_0),$$
which implies
\begin{equation}\label{eq:2.22}
\|u\|\lesssim \tau(\lambda)\left(\|v\|+\lambda^{1/2}\|f\|_0+\lambda^{1/2}\|\omega\|_0\right).
\end{equation} 
Now (\ref{eq:2.6}) follows from (\ref{eq:2.15}) and (\ref{eq:2.22}) in this case.
\eproof

Denote by $G$
the selfadjoint Dirichlet realization of the operator $-n(x)^{-1}\nabla c(x)\nabla$ on the Hilbert space 
$H=L^2(\Omega;n(x)dx)$. Then $G\ge 0$ and the spectrum of $G$ is discreet consisting of infinitely many non-negative eigenvalues. 
Given $\lambda\in\mathbb{C}^+$, set 
$$\theta(\lambda)={\rm dist}(\lambda,{\rm spec}\,\sqrt{G})=\min_{\mu_k\in {\rm spec}\,\sqrt{G}}|\lambda-\mu_k|\ge |{\rm Im}\,\lambda|.$$ 
In the next sections we will also need the following

\begin{Theorem} \label{2.6} 
Let $u\in H^2(\Omega)$ satisfy equation (\ref{eq:2.1}). 
Then the function $\omega=h\partial_\nu u|_\Gamma$ satisfies the estimate
\begin{equation}\label{eq:2.23}
\|\omega\|_0\lesssim \left(1+\theta(\lambda)^{-1}\right)\|v\|+\left(1+|\lambda|\theta(\lambda)^{-1}\right)\|f\|_{1,0}
\end{equation}
for $\theta(\lambda)>0$, $|{\rm Im}\,\lambda|\le C$, $|\lambda|\gg 1$, $C>0$ being any constant.
\end{Theorem}

{\it Proof.}  The theorem follows from the next two lemmas.

\begin{lemma} \label{2.7} 
We have the estimate
\begin{equation}\label{eq:2.24}
\|\omega\|_0\lesssim \|v\|+\|u\|+\|f\|_{1,0}.
\end{equation}
\end{lemma}

{\it Proof.} As in the proof of Lemma \ref{2.4}, it suffices to prove (\ref{eq:2.24}) locally near $\Gamma$. 
We keep the same notations. In the coordinates $(x_1,x')$ the principal symbol of the Euclidean Laplacian $-\Delta$ is equal to 
$\xi_1^2+r(x,\xi')$, where $(\xi_1,\xi')$ are the dual variables to $(x_1,x')$, and $r$ is a homogeneous polynomial of order two and satisfies $C_1|\xi'|^2\le r\le C_2|\xi'|^2$
with some constants $C_1,C_2>0$. 
Therefore, the principal symbol of the positive Laplace-Beltrami operator on $\Gamma$ is equal to
$r_0(x',\xi')=r(0,x',\xi')$. Note that $\Gamma$ can be considered as a Riemannian manifold without boundary with a 
Riemannian metric induced by the Euclidean one. Therefore, we can write the operator
$$P=-h^2n(x)^{-1}\nabla c(x)\nabla-1$$
in the coordinates $(x_1,x')$ as follows
$$P=\frac{c(x)}{n(x)}\left(\mathcal{D}_{x_1}^2+r(x,\mathcal{D}_{x'})\right)-1+h\mathcal{R}(x,\mathcal{D}_{x}),$$
where $\mathcal{D}_{x_1}=-ih\partial_{x_1}$, $\mathcal{D}_{x'}=-ih\partial_{x'}$, $\mathcal{D}_{x}=-ih\partial_{x}$ and 
$\mathcal{R}$ is a first-order differential operator. Introduce the function
$$F(x_1)=\left\|\mathcal{D}_{x_1}u^\flat\right\|_0^2
-\left\langle r(x_1,\cdot,\mathcal{D}_{x'})u^\flat,u^\flat\right\rangle_0+\left\langle n_\flat(x_1,\cdot)u^\flat,u^\flat\right\rangle_0,$$
where $n_\flat=c^{-1}n$. Clearly, 
\begin{equation}\label{eq:2.25}
{\rm Re}\,F(0)\ge\|\gamma\psi\omega\|_0^2-C\|f\|_{1,0}^2, \quad C>0.
\end{equation}
 On the other hand,
\begin{equation}\label{eq:2.26}
F(0)=-\int_0^{\epsilon}F'(x_1)dx_1
\end{equation} 
where $F'$ denotes the first derivative with respect to $x_1$. We have
$$F'(x_1)=-2{\rm Re}\,\left\langle (\mathcal{D}_{x_1}^2+r-n_\flat)u^\flat,
\partial_{x_1}u^\flat\right\rangle_0-\left\langle (r'-n_\flat')u^\flat,u^\flat\right\rangle_0$$
$$=2h^{-1}{\rm Im}\,\left\langle n_\flat(P-h\mathcal{R})u^\flat,\mathcal{D}_{x_1}u^\flat\right\rangle_0
-\left\langle (r'-n_\flat')u^\flat,u^\flat\right\rangle_0.$$
Hence
\begin{equation}\label{eq:2.27}
|F'(x_1)|\lesssim h^{-2}\|Pu^\flat\|_0^2+\|u^\flat\|_{1,0}^2+\|\mathcal{D}_{x_1}u^\flat\|_0^2.
\end{equation} 
By (\ref{eq:2.26}) and (\ref{eq:2.27}) we obtain
\begin{equation}\label{eq:2.28}
{\rm Re}\,F(0)\le \int_0^{\epsilon}|F'(x_1)|dx_1\lesssim  h^{-2}\|Pu^\flat\|^2+\|u^\flat\|_1^2.
\end{equation}
By (\ref{eq:2.25}) and (\ref{eq:2.28}),
$$
\left\|\gamma\psi\omega\right\|_0\lesssim h^{-1}\|Pu^\flat\|+\|u^\flat\|_1+\|f\|_{1,0}$$
$$\lesssim 
\lambda\|P((1-\phi)u)\|+\|(1-\phi)u\|_1+\|f\|_{1,0}$$
$$\lesssim 
\lambda\|Pu\|+\|u\|_1+\|f\|_{1,0},$$
which implies 
\begin{equation}\label{eq:2.29}
\|\omega\|_0\lesssim \|v\|+\|u\|_1+\|f\|_{1,0}.
\end{equation}
It is easy to see now that (\ref{eq:2.24}) follows from (\ref{eq:2.14}) and (\ref{eq:2.29}).
\eproof

\begin{lemma} \label{2.8} 
Under the assumptions of Theorem \ref{2.6} we have the estimate
\begin{equation}\label{eq:2.30}
\|u\|\lesssim \theta(\lambda)^{-1}\|v\|+\left(1+|\lambda|\theta(\lambda)^{-1}\right)\|f\|_{1,0}.
\end{equation}
\end{lemma}

{\it Proof.} Let $\zeta\in C^\infty(\overline\Omega)$ be supported in a small neighbourhood
of $\Gamma$ and $\zeta=1$ on $\Gamma$.
Clearly, if $u$ is a solution to the equation (\ref{eq:2.1}), then $u-\zeta f$
is a solution to (\ref{eq:2.1}) with $f$ replaced by $0$ and $v$ replaced by
$$v-\lambda^{-1}(\nabla c(x)\nabla+\lambda^2n(x))\zeta f.$$
Hence 
\begin{equation}\label{eq:2.31}
u-\zeta f=-\lambda(G-\lambda^2)^{-1}n^{-1}v+(G-\lambda^2)^{-1}n^{-1}(\nabla c(x)\nabla+\lambda^2n(x))\zeta f=: u_1+u_2.
\end{equation}
We have 
\begin{equation}\label{eq:2.32}
\|u_1\|\lesssim |\lambda|\left\|(G-\lambda^2)^{-1}\right\|_{H\to H}\|v\|$$ 
$$\lesssim |\lambda|\max_{\mu_k\in {\rm spec}\,\sqrt{G}}\left|\mu_k^2-\lambda^2\right|^{-1}\|v\|$$
$$\lesssim \max_{\mu_k\in {\rm spec}\,\sqrt{G}}\left|\mu_k-\lambda\right|^{-1}\|v\|
\lesssim \theta(\lambda)^{-1}\|v\|.
\end{equation}
To bound the norm of $u_2$, observe that given any $w\in L^2(\Omega)$, in view of (\ref{eq:2.19}), we have
$$
\langle u_2,w\rangle=\langle (\nabla c(x)\nabla+\lambda^2n(x))\zeta f,n^{-1}(G-\overline{\lambda}^2)^{-1}w\rangle$$  
$$=\langle c(x)\nabla(\zeta f),\nabla (n^{-1}(G-\overline{\lambda}^2)^{-1}w)\rangle
+\lambda^2\langle\zeta f,n^{-1}(G-\overline{\lambda}^2)^{-1}w\rangle.$$
Hence
\begin{equation}\label{eq:2.33}
|\langle u_2,w\rangle|\lesssim |\lambda|^2\|\zeta f\|_1\|(G-\overline{\lambda}^2)^{-1}w\|_1.
\end{equation}
On the other hand, since the function $W=(G-\overline{\lambda}^2)^{-1}w$ satisfies equation (\ref{eq:2.1})
with $f=0$, $\lambda$ replaced by $\overline{\lambda}$ and $v$ replaced by $\lambda^{-1}w$, the Green formula implies the estimate
\begin{equation}\label{eq:2.34}
\|W\|_1\lesssim \|W\|+|\lambda|^{-2}\|w\|.
\end{equation}
As in (\ref{eq:2.32}), we get
\begin{equation}\label{eq:2.35}
\|W\|\lesssim |\lambda|^{-1}\theta(\lambda)^{-1}\|w\|+|\lambda|^{-2}\|w\|,
\end{equation}
where we have used the fact that $\theta(\overline{\lambda})=\theta(\lambda)$. 
Since $\|\zeta f\|_1\lesssim \|f\|_{0,1}$, it follows from (\ref{eq:2.33}), (\ref{eq:2.34}) and (\ref{eq:2.35})
that
$$ |\langle u_2,w\rangle|\lesssim (1+|\lambda|\theta(\lambda)^{-1})\|f\|_{0,1}\|w\|,$$
which clearly implies 
\begin{equation}\label{eq:2.36}
\|u_2\|\lesssim (1+|\lambda|\theta(\lambda)^{-1})\|f\|_{0,1}.
\end{equation}
Now (\ref{eq:2.30}) follows from (\ref{eq:2.31}), (\ref{eq:2.32}) and (\ref{eq:2.36}).
\eproof

\section{Study of the interior Helmholtz equation with damping}

In this section we will obtain a priori estimates for the solution to the equation
\begin{equation}\label{eq:3.1}
\left\{
\begin{array}{l}
(\nabla c(x)\nabla+\lambda^2n(x)+i\lambda m(x))u=\lambda v\quad \mbox{in}\quad\Omega,\\
u=f\quad\mbox{on}\quad\Gamma,\\
u=0\quad\mbox{on}\quad\Gamma_0,
\end{array}
\right.
\end{equation}
where $\lambda\in\mathbb{C}$, ${\rm Re}\,\lambda\gg 1$,  
$c,n,m\in C^\infty(\overline\Omega)$, $c>0$, $n>0$, and the function $m$ does not change the sign. 
We keep the same notations as in the previous section. We define in the same way the operator $\widetilde G$ and suppose
that its resolvent satisfies the condition (\ref{eq:2.2}). 
Our goal in this section is to prove the following

\begin{Theorem} \label{3.1} 
Suppose that $m$ satisfies (\ref{eq:1.4}). Let $u\in H^2(\Omega)$ satisfy equation (\ref{eq:3.1}) and set 
 $\omega=h\partial_\nu u|_\Gamma$. Then there are constants $C,\lambda_0>0$ such that for all $\lambda\in\mathbb{C}$ such that 
 \begin{equation}\label{eq:3.2}
 |{\rm Im}\,\lambda|\le C\tau(|\lambda|)^{-1},\quad {\rm Re}\,\lambda\ge\lambda_0, 
 \end{equation}
 we have the estimate
 \begin{equation}\label{eq:3.3}
\|u\|\lesssim \tau(|\lambda|)\left(\|v\|+\|mu\|\right),
\end{equation}
 while for
 \begin{equation}\label{eq:3.4}
 |{\rm Im}\,\lambda|\le C\tau(|\lambda|)^{-2},\quad {\rm Re}\,\lambda\ge\lambda_0, 
 \end{equation}
 we have the estimates
\begin{equation}\label{eq:3.5}
\|\omega\|_0+\|u\|_1\lesssim \tau(|\lambda|)^2\left(\|v\|+\|f\|_0\right)+\|f\|_{1,0},
\end{equation}
\begin{equation}\label{eq:3.6}
\|u\|\lesssim \tau(|\lambda|)^2\|v\|+\tau(|\lambda|)|{\rm Im}\,\langle c\omega,f\rangle_0|^{1/2},
\end{equation}
\begin{equation}\label{eq:3.7}
\|u\|_1\lesssim \tau(|\lambda|)^2\|v\|+h^{1/2}\|f\|_{1,0}+\tau(|\lambda|)
|{\rm Im}\,\langle c\omega,f\rangle_0|^{1/2}.
\end{equation}
Moreover, under the condition
\begin{equation}\label{eq:3.8}
{\rm sign}(m)\,{\rm Im}\,\langle c\omega,f\rangle_0\le 0
\end{equation}
we have the better estimate
\begin{equation}\label{eq:3.9}
\|u\|\lesssim \tau(|\lambda|)^2\|v\|.
\end{equation}
\end{Theorem}

{\it Proof.}  The estimate (\ref{eq:3.3}) is a consequence of (\ref{eq:2.5}). Indeed, in view of assumption (\ref{eq:1.4})
we can find $\epsilon>0$ such that $|m|\ge C$ on $\Omega_\epsilon$ with some constant $C>0$. Applying (\ref{eq:2.5}) with
$v$ replaced by $v-imu$, we get
$$\|u\|\lesssim \tau(|\lambda|)\left(\|v-imu\|+\|u\|_{L^2(\Omega_\epsilon)}\right)
\lesssim \tau(|\lambda|)\left(\|v\|+\|mu\|\right).$$
As in the proof of Theorem 2.1, it is easy to see that it suffices to prove the estimates 
(\ref{eq:3.5}), (\ref{eq:3.6}), (\ref{eq:3.7}) and (\ref{eq:3.9}) for real $\lambda\gg 1$.
The Green formula in this case takes the form
\begin{equation}\label{eq:3.10}
\langle(\lambda^2n+i\lambda m)u-\lambda v,u\rangle=
\langle -\nabla c\nabla u,u\rangle=\int_\Omega c|\nabla u|^2+\lambda\langle c\omega,f\rangle_0.
\end{equation}
Taking the imaginary part of this identity we get
\begin{equation}\label{eq:3.11}
\langle mu,u\rangle={\rm Im}\langle v,u\rangle+{\rm Im}\langle c\omega,f\rangle_0,
\end{equation}
which implies
\begin{equation}\label{eq:3.12}
\||m|^{1/2}u\|^2\le \|v\|\|u\|+{\rm sign}(m)\,{\rm Im}\,\langle c\omega,f\rangle_0.
\end{equation}
Since $|m|\lesssim |m|^{1/2}$, we obtain from (\ref{eq:3.3}) and (\ref{eq:3.12}),
$$\|u\|^2\lesssim \tau(\lambda)^2\|v\|^2+\tau(\lambda)^2\||m|^{1/2}u\|^2$$ 
$$
\lesssim \tau(\lambda)^2\|v\|\|u\|+\tau(\lambda)^2\|v\|^2+\tau(\lambda)^2{\rm sign}(m)\,{\rm Im}\,\langle c\omega,f\rangle_0,$$
which implies
\begin{equation}\label{eq:3.13}
\|u\|^2\lesssim \tau(\lambda)^4\|v\|^2+\tau(\lambda)^2{\rm sign}(m)\,{\rm Im}\,\langle c\omega,f\rangle_0.
\end{equation}
Clearly, (\ref{eq:3.6}) and (\ref{eq:3.9}) follow from (\ref{eq:3.13}). 
Furthermore, taking the real part of (\ref{eq:3.10}) leads to the estimate
\begin{equation}\label{eq:3.14}
\|u\|_1\lesssim h\|v\|+\|u\|+h^{1/2}\|f\|_0^{1/2}\|\omega\|_0^{1/2}$$
$$\lesssim h\|v\|+\|u\|+h^{1/2}\|f\|_0+h^{1/2}\|\omega\|_0.
\end{equation}
On the other hand, applying (\ref{eq:2.24}) with
$v$ replaced by $v-imu$, we get 
\begin{equation}\label{eq:3.15}
\|\omega\|_0\lesssim \|v\|+\|u\|+\|f\|_{1,0}.
\end{equation}
By (\ref{eq:3.14}) and (\ref{eq:3.15}), 
\begin{equation}\label{eq:3.16}
\|u\|_1\lesssim h^{1/2}\|v\|+\|u\|+h^{1/2}\|f\|_{1,0}.
\end{equation}
Clearly, the estimate (\ref{eq:3.7}) follows from (\ref{eq:3.6}) and (\ref{eq:3.16}).
To prove (\ref{eq:3.5}) observe that from (\ref{eq:3.15}) and (\ref{eq:3.16}) together with (\ref{eq:3.6}) we obtain
$$\|\omega\|_0+\|u\|_1\lesssim \|v\|+\|u\|+\|f\|_{1,0}
\lesssim \tau(|\lambda|)^2\|v\|+\|f\|_{1,0}+\tau(|\lambda|)\|f\|_0^{1/2}\|\omega\|_0^{1/2}$$
$$\lesssim \tau(|\lambda|)^2\|v\|+\|f\|_{1,0}+\varepsilon^{-1}\tau(|\lambda|)^2\|f\|_0+\varepsilon\|\omega\|_0$$
for every $0<\varepsilon\le 1$. We now absorb the last term by taking $\varepsilon$ small enough, independent of $\lambda$.
\eproof

\section{A priori estimates for the boundary data}

Given any $k\in \mathbb{R}$, $S^k(\Gamma)$ will denote the set of all functions 
$a\in C^\infty(T^*\Gamma)$ satisfying
$$\left|\partial_x^\alpha\partial_{\xi}^\beta a\right|\le C_{\alpha,\beta}(|\xi|+1)^{k-|\beta|}$$
for all multi-indices $\alpha$ and $\beta$.
Given a symbol $a\in S^k(\Gamma)$, ${\rm Op}_h(a)$ will denote the
$h-\Psi$DO defined by
$$\left({\rm Op}_h(a)f\right)(x)=(2\pi h)^{1-d}\int_{T^*\Gamma}e^{i\langle x-y,\xi\rangle/h}
a(x,\xi)f(y)d\xi dy.$$
Let $\chi\in C^\infty(T^*\Gamma)$ be of compact support, independent of $\lambda$. We will derive from Theorem 3.1 the following

\begin{prop} \label{4.1} 
Suppose that $m$ satisfies (\ref{eq:1.4}). Let $u\in H^2(\Omega)$ satisfy equation (\ref{eq:3.1}) and set 
 $\omega=h\partial_\nu u|_\Gamma$. Then for all $\lambda\in\mathbb{C}$ satisfying (\ref{eq:3.4}) we have the estimate
\begin{equation}\label{eq:4.1}
\|f\|_{1,0}\lesssim |\lambda|^{1/2}\tau(|\lambda|)^2\|v\|+|\lambda|^{1/2}\tau(|\lambda|)|{\rm Im}\,\langle c\omega,f\rangle_0|^{1/2}
+\|{\rm Op}_h(1-\chi)f\|_{1,0}.
\end{equation}
\end{prop}

{\it Proof.} By (\ref{eq:2.16}) together with (\ref{eq:3.6}) and (\ref{eq:3.7}) we get
\begin{equation}\label{eq:4.2}
\|f\|_0\lesssim \varepsilon^{-1}|\lambda|^{1/2}\tau(|\lambda|)^2\|v\|+\varepsilon^{-1}|\lambda|^{1/2}
\tau(|\lambda|)|{\rm Im}\,\langle c\omega,f\rangle_0|^{1/2}
+\varepsilon\|f\|_{1,0}.
\end{equation}
On the other hand, we have 
\begin{equation}\label{eq:4.3}
\|f\|_{1,0}\lesssim \|{\rm Op}_h(\chi)f\|_{1,0}+\|{\rm Op}_h(1-\chi)f\|_{1,0}
\lesssim \|f\|_0+\|{\rm Op}_h(1-\chi)f\|_{1,0}.
\end{equation}
We now combine (\ref{eq:4.2}) and (\ref{eq:4.3}) and take $\varepsilon$ small enough, independent of $\lambda$, in order to absorb the term
$\varepsilon \|f\|_{1,0}$ 
in the right-hand side. This clearly leads to the estimate (\ref{eq:4.1}). 
\eproof

Let $\chi\in C^\infty(T^*\Gamma)$ be supported in the hyperbolic region, ${\cal H}$, of the boundary value problem (\ref{eq:3.1}), that is,
$$ {\rm supp}\,\chi\subset{\cal H}:=\{(x',\xi')\in T^*\Gamma:r_0(x',\xi')<n_{\flat,0}(x')\},$$
where $n_{\flat,0}=n_\flat|_\Gamma$, $n_\flat=n/c$. 
With this choice of the function $\chi$ we will now prove the following improved estimates. 

\begin{prop} \label{4.2} 
Under the assumptions of Proposition \ref{4.1}, we have the estimates
\begin{equation}\label{eq:4.4}
\|{\rm Op}_h(\chi)f\|_0+\|{\rm Op}_h(\chi)\omega\|_0$$
$$\lesssim
\tau(|\lambda|)^2\|v\|+\tau(|\lambda|)|{\rm Im}\,
\langle c\omega,f\rangle_0|^{1/2}+h^{1/2}\|f\|_{1,0},
\end{equation}
\begin{equation}\label{eq:4.5}
\|f\|_{1,0}\lesssim
\tau(|\lambda|)^2\|v\|+\tau(|\lambda|)|{\rm Im}\,\langle c\omega,f\rangle_0|^{1/2}
+\|{\rm Op}_h(1-\chi)f\|_{1,0}.
\end{equation}
\end{prop}

{\it Proof.} By (\ref{eq:4.3}) we have 
\begin{equation}\label{eq:4.6}
\|f\|_{1,0}\lesssim \|{\rm Op}_h(\chi)f\|_0+\|{\rm Op}_h(1-\chi)f\|_{1,0}.
\end{equation}
Clearly, the estimate (\ref{eq:4.5}) follows from (\ref{eq:4.4}) and (\ref{eq:4.6}). On the other hand, 
it is easy to see that (\ref{eq:4.4}) follows from (\ref{eq:3.7})  
and the following 

\begin{lemma} \label{4.3} 
We have the estimate
\begin{equation}\label{eq:4.7}
\|{\rm Op}_h(\chi)f\|_0+\|{\rm Op}_h(\chi)\omega\|_0\lesssim\|v\|+\|u\|_1+h\|f\|_0.
\end{equation}
\end{lemma}

{\it Proof.} It suffices to prove the estimate (\ref{eq:4.7}) for real $\lambda\gg 1$. Indeed, this would imply (\ref{eq:4.7})
for complex $\lambda$ such that $|{\rm Im}\,\lambda|\le C$, $C>0$ being any constant. 
Let the function $u^\flat$ and the operator $P$ be as in the proof of Lemma \ref{2.7}. Set 
$$u^\sharp={\rm Op}_h(\chi)u^\flat={\rm Op}_h(\chi)\psi(1-\phi)u.$$ 
Clearly,
\begin{equation}\label{eq:4.8}
\|Pu^\sharp\|\lesssim \|P((1-\phi)u)\|+\|[P,{\rm Op}_h(\chi)\psi](1-\phi)u\|$$ 
$$\lesssim \|P((1-\phi)u)\|+h\|(1-\phi)u\|_1$$
$$\lesssim \|Pu\|+h\|u\|_1\lesssim h\|v\|+h\|u\|_1.
\end{equation}
We define the function $F(x_1)$ as in the proof of Lemma \ref{2.7} replacing $u^\flat$ by $u^\sharp$. Observe now that the choice of
$\chi$ guarantees that
$$n_\flat(0,x')-r_0(x',\xi')\ge C>0$$
on supp$\,\chi$. Therefore, by G{\aa}rding's inequality we have 
\begin{equation}\label{eq:4.9}
{\rm Re}\,\left\langle (n_\flat(0,\cdot)- r_0(\cdot,{\cal D}_{x'})){\rm Op}_h(\chi)\psi_0f,{\rm Op}_h(\chi)\psi_0f\right\rangle_0
\ge C_1\|{\rm Op}_h(\chi)\psi_0f\|_0^2,
\end{equation}
where $C_1>0$ is some constant and $\psi_0=\psi|_{x_1=0}$. 
Since
$${\cal D}_{x_1}u^\sharp|_{x_1=0}=-i{\rm Op}_h(\chi)\psi_0\omega-ih{\rm Op}_h(\chi)\psi_1f,$$
where $\psi_1=\partial_{x_1}\psi|_{x_1=0}$, we deduce from (\ref{eq:4.9}),
\begin{equation}\label{eq:4.10}
{\rm Re}\,F(0)\ge C_1\|{\rm Op}_h(\chi)\psi_0f\|_0^2+\|{\rm Op}_h(\chi)\psi_0\omega\|_0^2-O(h^2)\|f\|_0^2$$ 
$$\ge C_1\|\psi_0{\rm Op}_h(\chi)f\|_0^2+\|\psi_0{\rm Op}_h(\chi)\omega\|_0^2-O(h^2)\|f\|_0^2.
\end{equation}
On the other hand, the upper bound (\ref{eq:2.28}) still holds with $u^\flat$ replaced by $u^\sharp$.
This fact together with (\ref{eq:4.8}) and (\ref{eq:4.10}) imply (\ref{eq:4.7}).
\eproof

Given a parameter $0<\varepsilon\ll 1$, independent of $\lambda$, choose a function  $\chi_\varepsilon\in C_0^\infty(T^*\Gamma)$
such that $\chi_\varepsilon=1$ in the region $\{|r_0/n_{\flat,0}-1|\le\varepsilon\}$ and 
$\chi_\varepsilon=0$ in $T^*\Gamma\setminus\{|r_0/n_{\flat,0}-1|\le 2\varepsilon\}$.
We have the following

\begin{prop} \label{4.4} 
Under the assumptions of Proposition \ref{4.1}, we have the estimate
\begin{equation}\label{eq:4.11}
\|{\rm Op}_h(\chi_\varepsilon)\omega\|_0\lesssim 
\tau(|\lambda|)^2\|v\|+\tau(|\lambda|)|{\rm Im}\,\langle c\omega,f\rangle_0|^{1/2}+(\varepsilon+h^{1/2})\|f\|_{1,0}.
\end{equation}
\end{prop}

{\it Proof.} It is easy to see that the proposition follows from (\ref{eq:3.7}) and the following 

\begin{lemma} \label{4.5} 
We have the estimate
\begin{equation}\label{eq:4.12}
\|{\rm Op}_h(\chi_\varepsilon)\omega\|_0\lesssim\|v\|+\|u\|_1+(\varepsilon+h)\|f\|_0.
\end{equation}
\end{lemma}

{\it Proof.} Again, it suffices to prove the lemma for real $\lambda\gg 1$.
We will proceed in the same way as in the proof of Lemma \ref{4.3} with $\chi$ replaced by $\chi_\varepsilon$ making the following modification. Since in this case the function $\chi_\varepsilon$ is no longer supported in the hyperbolic region, we do not have
the G{\aa}rding inequality (\ref{eq:4.9}) fulfilled anymore. Instead, since $n_\flat(0,x')-r_0(x',\xi')=O(\varepsilon)$ on 
supp$\,\chi_\varepsilon$, we have the bound
$$\left\|(n_\flat(0,\cdot)- r_0(\cdot,\mathcal{D}_{x'})){\rm Op}_h(\chi_\varepsilon)\psi_0f\right\|_0
\lesssim (\varepsilon+h)\|\psi_0f\|_0.$$
Therefore in this case the function ${\rm Re}\,F(0)$ is lower bounded as follows
\begin{equation}\label{eq:4.13}
{\rm Re}\,F(0)\ge \|{\rm Op}_h(\chi_\varepsilon)\psi_0\omega\|_0^2-O((\varepsilon+h)^2)\|f\|_0^2$$ 
$$\ge \|\psi_0{\rm Op}_h(\chi_\varepsilon)\omega\|_0^2-O((\varepsilon+h)^2)\|f\|_0^2.
\end{equation}
The estimate (\ref{eq:4.12}) follows from (\ref{eq:4.13}) and the upper bound (\ref{eq:2.28}) adapted to this case.
\eproof

\section{The Dirichlet-to-Neumann map}

Let $u$ solve equation (\ref{eq:3.1}) with $v\equiv 0$ and define the Dirichlet-to-Neumann map
$$\mathcal{N}(\lambda,m):H^1(\Gamma)\to L^2(\Gamma)$$
by
$$\mathcal{N}(\lambda,m)f:=h\partial_\nu u|_\Gamma.$$
When $m\equiv 0$ we will denote the Dirichlet-to-Neumann map by $\mathcal{N}(\lambda)$. 
Clearly, under the conditions of Theorem \ref{3.1}, by (\ref{eq:3.5}) we have that 
the Dirichlet-to-Neumann map in this case satisfies the estimate 
\begin{equation}\label{eq:5.1}
\|\mathcal{N}(\lambda,m)f\|_0\lesssim \tau(|\lambda|)^2\|f\|_{1,0}
\end{equation}
for $\lambda$ satisfying (\ref{eq:3.4}). 
On the other hand, when $m\equiv 0$, by (\ref{eq:2.23}) we have
the estimate
\begin{equation}\label{eq:5.2}
\|\mathcal{N}(\lambda)f\|_0\lesssim\left(1+|\lambda|\theta(\lambda)^{-1}\right)\|f\|_{1,0},
\end{equation}
for $\theta(\lambda)>0$, $|{\rm Im}\,\lambda|\le C$, $|\lambda|\gg 1$, $C>0$ being any constant.
Let $\chi,\eta\in C^\infty(T^*\Gamma)$ be compactly supported functions such that $\eta=1$ on supp$\,\chi$. 
 In Section 7 we will need the following

\begin{lemma} \label{5.1} 
Under the assumptions of Theorem \ref{3.1}, we have the estimates
\begin{equation}\label{eq:5.3}
\|[\mathcal{N}(\lambda,m),{\rm Op}_h(\chi)]f\|_0\lesssim h^{1/2}\tau(|\lambda|)^2\|f\|_{1,0}
+\tau(|\lambda|)^3|{\rm Im}\,\langle c\mathcal{N}(\lambda,m)f,f\rangle_0|^{1/2},
\end{equation}
\begin{equation}\label{eq:5.4}
\|{\rm Op}_h(1-\eta)\mathcal{N}(\lambda,m){\rm Op}_h(\chi)f\|_0\lesssim h^{1/2}\tau(|\lambda|)^2\|f\|_{1,0}
+\tau(|\lambda|)^3|{\rm Im}\,\langle c\mathcal{N}(\lambda,m)f,f\rangle_0|^{1/2}.
\end{equation}
\end{lemma}

{\it Proof.} Let $u$ solve equation (\ref{eq:3.1}) with $v\equiv 0$. Then the function
$\widetilde u={\rm Op}_h(\chi)(1-\phi)u$ solves equation (\ref{eq:3.1}) with $v$ and $f$ replaced by $\widetilde v$
and $\widetilde f$, respectively, where
$$\widetilde v=\lambda^{-1}[\nabla c(x)\nabla+\lambda^2n(x)+i\lambda m(x),{\rm Op}_h(\chi)(1-\phi)]u,$$
$$\widetilde f={\rm Op}_h(\chi)f.$$
Let $w$ solve equation (\ref{eq:3.1}) with $v\equiv 0$ and $f$ replaced by $\widetilde f$.
Then the function $\widetilde u-w$ solves equation (\ref{eq:3.1}) with $v$ replaced by $\widetilde v$ and $f=0$.
Moreover, we have
$$h\partial_\nu(\widetilde u-w)|_\Gamma=-[\mathcal{N}(\lambda,m),{\rm Op}_h(\chi)]f.$$
Therefore, the estimates (\ref{eq:3.5}) and (\ref{eq:3.7}) lead to 
$$\|[\mathcal{N}(\lambda,m),{\rm Op}_h(\chi)]f\|_0\lesssim \tau(|\lambda|)^2\|\widetilde v\|\lesssim \tau(|\lambda|)^2\|u\|_1$$ 
$$\lesssim h^{1/2}\tau(|\lambda|)^2\|f\|_{1,0}
+\tau(|\lambda|)^3|{\rm Im}\,\langle c\mathcal{N}(\lambda,m)f,f\rangle_0|^{1/2}.$$
To prove (\ref{eq:5.4}) we will use that
$${\rm Op}_h(1-\eta){\rm Op}_h(\chi)=O(h^\infty):L^2(\Gamma)\to L^2(\Gamma).$$
Thus, in view of (\ref{eq:5.1}), we get
$$\|{\rm Op}_h(1-\eta)\mathcal{N}(\lambda,m){\rm Op}_h(\chi)f\|_0\lesssim 
\|{\rm Op}_h(1-\eta)[\mathcal{N}(\lambda,m),{\rm Op}_h(\chi)]f\|_0$$
$$+\|{\rm Op}_h(1-\eta){\rm Op}_h(\chi)\mathcal{N}(\lambda,m)f\|_0$$
$$\lesssim \|[\mathcal{N}(\lambda,m),{\rm Op}_h(\chi)]f\|_0+h^\infty\|\mathcal{N}(\lambda,m)f\|_0$$
$$\lesssim \|[\mathcal{N}(\lambda,m),{\rm Op}_h(\chi)]f\|_0+h^\infty\tau(|\lambda|)^2\|f\|_{1,0}.$$
We now obtain (\ref{eq:5.4}) from (\ref{eq:5.3}).
\eproof

Denote by $\mathcal{N}(\lambda,m)^*$ the adjoint of $\mathcal{N}(\lambda,m)$ with respect to the scalar product
$\langle\cdot,\cdot\rangle_0$ in $L^2(\Gamma)$. 
In Section 6 we will need the following

\begin{lemma} \label{5.2} The adjoint operator of the Dirichlet-to-Neumann map satisfies the identity
\begin{equation}\label{eq:5.5}
\mathcal{N}(\lambda,m)^*c_0=c_0\mathcal{N}(\overline\lambda,-m),
\end{equation}
where $c_0=c|_\Gamma$. 
\end{lemma}

{\it Proof.} Given any $f_1,f_2\in L^2(\Gamma)$, let $u_1$ be the solution of equation (\ref{eq:3.1})
with $v\equiv 0$ and $f$ replaced by $f_1$, and let $u_2$ be the solution of equation (\ref{eq:3.1})
with $v\equiv 0$, $f$ replaced by $f_2$ and $\lambda$ replaced by $\overline\lambda$. By the Green formula we have 
$$0=-\langle\nabla c\nabla u_1,u_2\rangle+\langle u_1,\nabla c\nabla u_2\rangle
=\langle c_0\partial_\nu u_1|_\Gamma,f_2\rangle_0-\langle f_1,c_0\partial_\nu u_2|_\Gamma\rangle_0$$
 $$=h^{-1}\langle c_0\mathcal{N}(\lambda,m)f_1,f_2\rangle_0-h^{-1}\langle f_1,c_0\mathcal{N}(\overline\lambda,-m)f_2\rangle_0,$$
which clearly implies (\ref{eq:5.5}).
\eproof

\section{Parametrix of the Dirichlet-to-Neumann map in the elliptic region revisited}

Let $\eta\in C^\infty(T^*\Gamma)$ be such that $1-\eta$ is supported in the elliptic region, ${\cal E}$, of the boundary value problem (\ref{eq:3.1}), that is,
$${\rm supp}(1-\eta)\subset{\cal E}:=\{(x',\xi')\in T^*\Gamma:r_0(x',\xi')>n_{\flat,0}(x')\}.$$
For $(x',\xi')\in {\cal E}$ set 
$$\rho(x',\xi',z)=\sqrt{r_0(x',\xi')-zn_{\flat,0}(x')},\quad {\rm Re}\,\rho>0,$$
where $z=(h\lambda)^2=1-(h{\rm Im}\,\lambda)^2+2ih{\rm Im}\,\lambda$. On ${\rm supp}(1-\eta)$ we have the lower bound
\begin{equation}\label{eq:6.1}
{\rm Re}\,\rho\ge C\langle\xi'\rangle, \quad C>0.
\end{equation} 
The main result in this section is the following

\begin{Theorem} \label{6.1} Suppose that the function $\tau$ satisfies the bound 
\begin{equation}\label{eq:6.2}
\tau(\lambda)\le \lambda^q
\end{equation}
with some constant $q\ge 0$. Then, under the conditions of Theorem \ref{3.1} 
we have the estimate
\begin{equation}\label{eq:6.3}
\left\|\mathcal{N}(\lambda,m){\rm Op}_h(1-\eta)f+{\rm Op}_h(\rho(1-\eta))f\right\|_0\lesssim h\|f\|_0 
\end{equation}
for $\lambda$ satisfying (\ref{eq:3.4}). When $m\equiv 0$ the estimate (\ref{eq:6.3}) still holds without assuming
(\ref{eq:6.2}) for 
$\lambda\in\mathbb{C}^+\setminus{\cal L}_N^+(C_N)$,
$|{\rm Im}\,\lambda|\le C$, $C>0$ being any constant, while $C_N>0$ is a suitable constant depending on $N$.
 \end{Theorem}

{\it Proof.} The theorem follows from the parametrix construction carried out in \cite{kn:V1}. 
In what follows we will recall it (see also Section 5 of \cite{kn:V2}).
In fact, in \cite{kn:V1} the case $m\equiv 0$ is considered, but it is easy to see that the presence of the function
$m$ does not change anything. Indeed, the eikonal equation does not depend on $m$ and only the transport equations do.
Note also that it suffices to build the parametrix locally and then sum up all pieces.

  Let $(x_1,x')\in {\cal V}^+$ be the local normal geodesic coordinates near the boundary. 
Take a function $\chi\in C^\infty(T^*\Gamma)$, $0\le\chi\le 1$, such that $\pi_{x'}({\rm supp}\,\chi)\subset {\cal V}^0$, where 
$\pi_{x'}:T^*\Gamma\to\Gamma$ denotes the projection $(x',\xi')\to x'$. Moreover, we require that 
 $\chi\in S^0(\Gamma)$ with ${\rm supp}\,\chi\subset{\rm supp}(1-\eta)$. 
We will be looking for a parametrix of the solution to 
equation (\ref{eq:3.1}) (with $v\equiv 0$) in the form
$$\widetilde u=\phi_0(x_1)(2\pi h)^{-d+1}\int\int e^{\frac{i}{h}(\langle y',\xi'\rangle+\varphi(x,\xi',z))}a(x,\xi',z,h)f(y')d\xi'dy',$$
where $\phi_0\in C_0^\infty(\mathbb{R})$,
$\phi_0(t)=1$ for $|t|\le \delta/2$, $\phi_0(t)=0$ for $|t|\ge\delta$. Here $0<\delta\ll 1$ is a small parameter independent of $\lambda$.   
 We require that $\widetilde u$ satisfies the boundary condition
$\widetilde u={\rm Op}_h(\chi)f$ on $x_1=0$. The phase $\varphi$ and the amplitude $a$ are choosen in such a way that the function
$\widetilde u$ satisfies equation (\ref{eq:3.1}) mod ${\cal O}(h^M)$, where $M\ge 1$ is an arbitrary integer. 
The phase function satisfies 
$$\varphi|_{x_1=0}=-\langle x',\xi'\rangle$$
as well as the eikonal equation
\begin{equation}\label{eq:6.4}
(\partial_{x_1}\varphi)^2+r(x,\nabla_{x'}\varphi)-z\widetilde n(x)=x_1^M\Psi_M,
\end{equation}
where the function $\left|\Psi_M\right|$ is bounded as $x_1\to 0$. It is shown in Section 4 of \cite{kn:V1} that 
(\ref{eq:6.4}) has a solution of the form
$$\varphi=\sum_{j=0}^Mx_1^j\varphi_j,$$
where the functions $\varphi_j$ do not depend on $x_1$, $\varphi_0=-\langle x',\xi'\rangle$,
$\varphi_1=i\rho$. It follows from (\ref{eq:6.1}) that
\begin{equation}\label{eq:6.5}
{\rm Im}\,\varphi\ge Cx_1\langle\xi'\rangle/2,
\end{equation} 
for $0\le x_1\le\delta$, provided $\delta$ is taken small enough. 
The amplitude is of the form
$$a=\sum_{j=0}^Mh^ja_j,$$
where the functions $a_j$ do not depend on $h$, $a_0|_{x_1=0}=\chi$. 
Then all functions $a_j$ can be determined from the transport equations
and we have $a_j\in S^{-j}(\Gamma)$ uniformly in $x_1$ and $z$ (see Section 4 of \cite{kn:V1}). Clearly, we have
\begin{equation}\label{eq:6.6}
h\partial_\nu\widetilde u|_{x_1=0}={\rm Op}_h(b_M)f,
\end{equation}
where
$$b_M=ia\frac{\partial\varphi}{\partial x_1}|_{x_1=0}+h\frac{\partial a}{\partial x_1}|_{x_1=0}=
-\chi\rho+h\sum_{j=0}^Mh^j\frac{\partial a_j}{\partial x_1}|_{x_1=0}.$$
Hence $h^{-1}(b_M+\chi\rho)\in S^0(\Gamma)$ uniformly in $h$. This implies
\begin{equation}\label{eq:6.7}
{\rm Op}_h(b_M+\chi\rho)=O(h):L^2(\Gamma)\to L^2(\Gamma).
\end{equation}
On the other hand, the function
$$\widetilde v=(\nabla c(x)\nabla+\lambda^2n(x)+i\lambda m(x))\widetilde u$$
is of the form
$$\widetilde v=(2\pi h)^{-d+1}\int\int e^{\frac{i}{h}(\langle y',\xi'\rangle+\varphi(x,\xi',z))}V_M(x,\xi',z,h)f(y')d\xi'dy',$$
where $V_M=V_M^{(1)}+\phi_0(x_1)V_M^{(2)}$, 
$$V_M^{(1)}=[\nabla c(x)\nabla,\phi_0(x_1)]a,$$
$$V_M^{(2)}=e^{-i\varphi/h}(\nabla c(x)\nabla+\lambda^2n(x)+i\lambda m(x))e^{i\varphi/h}a.$$
As shown in Section 4 of \cite{kn:V1}, the functions $a_j$ can be choosen in such a way that the function $V_M^{(2)}$ is of the form
\begin{equation}\label{eq:6.8}
V_M^{(2)}=x_1^MA_M+h^MB_M,
\end{equation}
where $A_M$ and $B_M$ are smooth functions. More precisely, since $\chi$ is supported in the elliptic region, we have 
$A_M\in S^2(\Gamma)$, $B_M\in S^{1-M}(\Gamma)$ uniformly in $h$, $z$ and $0<x_1\le\delta$
(see Proposition 3.4 of \cite{kn:V1}). Note that in view of (\ref{eq:6.5}) we have the bound
$$\left|x_1^Me^{i\varphi/h}\right|\lesssim h^M\langle\xi'\rangle^{-M}.$$
Thus we get that the function $V_M$ satisfies the bound
\begin{equation}\label{eq:6.9}
\left|V_Me^{i\varphi/h}\right|\lesssim h^M\langle\xi'\rangle^{-M+1}.
\end{equation}
By (\ref{eq:6.9}) we obtain the estimate
\begin{equation}\label{eq:6.10}
\|\widetilde v\|\le \widetilde C_Mh^{M/2}\|f\|_0,
\end{equation}
provided $M$ is taken big enough. Let $u$ solve equation (\ref{eq:3.1}) with $v\equiv 0$ and $f$ replaced by ${\rm Op}_h(\chi)f$. 
Then the function
$\widetilde u-u$ solves equation (\ref{eq:3.1}) with $v$ replaced by $\widetilde v$ and $f=0$. Therefore, under the conditions of
Theorem 3.1, by (\ref{eq:3.5}) together with (\ref{eq:6.2}), (\ref{eq:6.6}) and (\ref{eq:6.10}) we get
\begin{equation}\label{eq:6.11}
\left\|\mathcal{N}(\lambda,m){\rm Op}_h(\chi)f-{\rm Op}_h(b_M)f\right\|_0=
\|h\partial_{x_1}(\widetilde u-u)|_{x_1=0}\|\lesssim\tau(|\lambda|)^2\|\widetilde v\|\lesssim h^{M/2-2q}\|f\|_0.
\end{equation}
Taking $M\ge 4q+2$, by (\ref{eq:6.7}) and (\ref{eq:6.11}), we get
\begin{equation}\label{eq:6.12}
\left\|\mathcal{N}(\lambda,m){\rm Op}_h(\chi)f+{\rm Op}_h(\chi\rho)f\right\|_0\lesssim h\|f\|_0,
\end{equation}
which implies (\ref{eq:6.3}) in this case since $1-\eta$ can be written as a finite sum of functions $\chi$ for which
(\ref{eq:6.12}) holds. Consider now the case when $m\equiv 0$. We proceed similarly with the difference that 
we use the estimate (\ref{eq:2.23}) instead of (\ref{eq:3.5}). For $\lambda\in\mathbb{C}^+\setminus\mathcal{L}_N^+(C_N)$, 
$|{\rm Im}\,\lambda|\le C$, we obtain
\begin{equation}\label{eq:6.13}
\left\|\mathcal{N}(\lambda){\rm Op}_h(\chi)f-{\rm Op}_h(b_M)f\right\|_0\lesssim\left(1+\theta(\lambda)^{-1}\right)\|\widetilde v\|
\le \widetilde C_M^\sharp C_N^{-1}h^{M/2-N}\|f\|_0\le h\|f\|_0,
\end{equation}
provided we take $M=2N+1$ and $C_N=\widetilde C_M^\sharp$. Thus we conclude that the estimate (\ref{eq:6.12}) (and hence (\ref{eq:6.3}))
 still holds in this case as long as
$\lambda\in\mathbb{C}^+\setminus\mathcal{L}_N^+(C_N)$, $|{\rm Im}\,\lambda|\le C$.
\eproof

Let $\chi,\eta\in C^\infty(T^*\Gamma)$ be compactly supported functions such that $\eta=1$ on supp$\,\chi$
and supp$(1-\eta)\subset{\cal E}$. We will use Theorem \ref{6.1} to prove the following

\begin{lemma} \label{6.2} 
Under the conditions of Theorem \ref{6.1} we have the estimates
\begin{equation}\label{eq:6.14}
\|{\rm Op}_h(\chi)\mathcal{N}(\lambda,m){\rm Op}_h(1-\eta)f\|_0\lesssim h\|f\|_0,
\end{equation}
\begin{equation}\label{eq:6.15}
\|{\rm Op}_h(1-\eta)\mathcal{N}(\lambda,m){\rm Op}_h(\chi)f\|_0\lesssim h\|f\|_0.
\end{equation}
When $m\equiv 0$ the estimates (\ref{eq:6.14}) and (\ref{eq:6.15}) still hold for 
$\lambda\in\mathbb{C}^+\setminus{\cal L}_N^+(C_N)$, $|{\rm Im}\,\lambda|\le C$.
\end{lemma}

{\it Proof.} Since
$${\rm Op}_h(\chi){\rm Op}_h(\rho(1-\eta))=O(h^\infty):L^2(\Gamma)\to L^2(\Gamma),$$
the estimate (\ref{eq:6.14}) follows from (\ref{eq:6.3}). In view of Lemma \ref{5.2} the adjoint of the operator
$$\mathcal{A}:={\rm Op}_h(1-\eta)\mathcal{N}(\lambda,m){\rm Op}_h(\chi)$$
is
$$\mathcal{A}^*={\rm Op}_h(\chi)^*c_0\mathcal{N}(\overline\lambda,-m)c_0^{-1}{\rm Op}_h(1-\eta)^*.$$
Choose compactly supported functions $\chi_1,\eta_1\in C^\infty(T^*\Gamma)$ such that $\eta_1=1$ on supp$\,\chi_1$, 
 supp$(1-\eta_1)\subset{\cal E}$, $\chi_1=1$ on supp$\,\chi$ and $\eta=1$ on supp$\,\eta_1$.
 The standard $h-\Psi$DO calculus give
 $${\rm Op}_h(\chi)^*c_0{\rm Op}_h(1-\chi_1)=O(h^\infty):H^{-1}(\Gamma)\to L^2(\Gamma),$$
 $${\rm Op}_h(\eta_1)c_0^{-1}{\rm Op}_h(1-\eta)^*=O(h^\infty):L^2(\Gamma)\to H^1(\Gamma).$$
 We now apply the estimate (\ref{eq:6.14}) with $\chi$, $\eta$, $\lambda$ replaced by 
 $\chi_1$, $\eta_1$, $\overline\lambda$, respectively. We will also use (\ref{eq:6.3}) with $\eta$, $\lambda$, $m$ and $\rho$ replaced by
 $\eta_1$, $\overline\lambda$, $-m$ and $\overline\rho$, respectively. Note that $\overline\rho(1-\eta_1)\in S^1(\Gamma)$. 
 Thus, in view of (\ref{eq:5.1}) and (\ref{eq:6.2}), we get
 $$\|\mathcal{A}^*f\|_0\lesssim \|{\rm Op}_h(\chi_1)\mathcal{N}(\overline\lambda,-m)
 {\rm Op}_h(1-\eta_1)c_0^{-1}{\rm Op}_h(1-\eta)^*f\|_0$$
 $$+\|\mathcal{N}(\overline\lambda,-m){\rm Op}_h(\eta_1)c_0^{-1}{\rm Op}_h(1-\eta)^*f\|_0$$
 $$+\|{\rm Op}_h(\chi)^*c_0{\rm Op}_h(1-\chi_1)\mathcal{N}(\overline\lambda,-m)c_0^{-1}{\rm Op}_h(1-\eta)^*f\|_0$$
 $$\lesssim h\|c_0^{-1}{\rm Op}_h(1-\eta)^*f\|_0$$
 $$+|\lambda|^{2q}\|{\rm Op}_h(\eta_1)c_0^{-1}{\rm Op}_h(1-\eta)^*f\|_{1,0}$$
 $$+h^\infty\|\mathcal{N}(\overline\lambda,-m){\rm Op}_h(\eta_1)c_0^{-1}{\rm Op}_h(1-\eta)^*f\|_0$$
 $$+h^\infty\|\mathcal{N}(\overline\lambda,-m){\rm Op}_h(1-\eta_1)c_0^{-1}{\rm Op}_h(1-\eta)^*f\|_{-1,0}$$
 $$\lesssim h\|f\|_0+h^\infty\|f\|_0$$
 $$+h^\infty\|{\rm Op}_h(\overline\rho(1-\eta_1))c_0^{-1}{\rm Op}_h(1-\eta)^*f\|_{-1,0}$$
 $$+h^\infty\|(\mathcal{N}(\overline\lambda,-m){\rm Op}_h(1-\eta_1)+{\rm Op}_h(\overline\rho(1-\eta_1)))c_0^{-1}{\rm Op}_h(1-\eta)^*f\|_0$$
 $$\lesssim h\|f\|_0+h^\infty\|f\|_0\lesssim h\|f\|_0,$$
 where $\|\cdot\|_{-1,0}$ denotes the semiclassical norm in $H^{-1}(\Gamma)$. 
 In other words,
 $$\mathcal{A}^*=O(h):L^2(\Gamma)\to L^2(\Gamma),$$
 and hence so is the operator $\mathcal{A}$. Clearly, the same analysis still holds when
 $m\equiv 0$, using (\ref{eq:5.2}) instead of (\ref{eq:5.1}).
 \eproof
  
\section{Eigenvalue-free regions}

Let $(u_1,u_2)$ be the solution to equation (\ref{eq:1.9}) and set $f=u_1|_\Gamma=u_2|_\Gamma$. 
Then we can express the restrictions of the normal derivative of $u_1$ and $u_2$ in terms of 
the corresponding Dirichlet-to-Neumann maps, that is,
$$h\partial_\nu u_1|_\Gamma=\mathcal{N}_1(\lambda,m)f,\quad h\partial_\nu u_2|_\Gamma=\mathcal{N}_2(\lambda)f.$$
Therefore, $\lambda$ is a transmission eigenvalue if $T(\lambda)f\equiv 0$, where
$$T(\lambda)=c_1\mathcal{N}_1(\lambda,m)-c_2\mathcal{N}_2(\lambda).$$
We have to show that, for $\lambda$ belonging to the eigenvalue-free regions of Theorems \ref{1.1} and \ref{1.2}, if $T(\lambda)f\equiv 0$,
then $f\equiv 0$. Without loss of generality we may suppose that $\lambda\in\mathbb{C}^+$ (see Remark \ref{1.3}). 
We will first prove the following

\begin{lemma} \label{7.1}
If $\lambda\in\mathbb{C}^+$ satisfies the condition
\begin{equation}\label{eq:7.1}
{\rm sign}(m)\,{\rm Im}\,\lambda\le 0,
\end{equation}
then 
\begin{equation}\label{eq:7.2}
{\rm sign}(m)\,{\rm Im}\,\langle c_1\mathcal{N}_1(\lambda,m)f,f\rangle_0\le 0.
\end{equation}
Moreover, 
there exist constants $C,\lambda_0>0$ such that for all $\lambda\in\mathbb{C}^+$ satisfying 
\begin{equation}\label{eq:7.3}
 |{\rm Im}\,\lambda|\le C\tau_1(|\lambda|)^{-2},\quad {\rm Re}\,\lambda\ge\lambda_0,
 \end{equation}
 \begin{equation}\label{eq:7.4}
 |{\rm Im}\,\lambda|\le C\tau_2(|\lambda|)^{-1},\quad {\rm Re}\,\lambda\ge\lambda_0,
 \end{equation}
we have the estimate
\begin{equation}\label{eq:7.5}
|{\rm Im}\,\langle c_1\mathcal{N}_1(\lambda,m)f,f\rangle_0|\lesssim |{\rm Im}\,\lambda||\lambda|^{\ell_2}
\tau_1(|\lambda|)^4\tau_2(|\lambda|)^2\|f\|_{1,0}^2.
\end{equation}
\end{lemma}

{\it Proof.} The Green formula applied to the second equation in (\ref{eq:1.9}) gives the identity
$$2{\rm Im}\,\lambda\langle n_2u_2,u_2\rangle={\rm Im}\,\langle c_2\mathcal{N}_2(\lambda)f,f\rangle_0={\rm Im}\,\langle c_1
\mathcal{N}_1(\lambda,m)f,f\rangle_0.$$
Obviously, (\ref{eq:7.1}) implies (\ref{eq:7.2}). 
On the other hand, by (\ref{eq:2.6}) and (\ref{eq:5.1}) we have
$$\|u_2\|\lesssim |\lambda|^{\ell_2/2}\tau_2(|\lambda|)\left(\|f\|_0+\|\mathcal{N}_2(\lambda)f\|_0\right)$$
 $$\lesssim |\lambda|^{\ell_2/2}\tau_2(|\lambda|)\left(\|f\|_0+\|\mathcal{N}_1(\lambda,m)f\|_0\right)$$
 $$\lesssim |\lambda|^{\ell_2/2}\tau_2(|\lambda|)\tau_1(|\lambda|)^2\|f\|_{1,0},$$
which clearly implies (\ref{eq:7.5}).
\eproof

To prove Theorem \ref{1.1} we apply Theorem \ref{3.1} with $u=u_1$ and $v=0$. Note that the conditions (\ref{eq:3.8})
and (\ref{eq:7.2}) are equivalent. Therefore, the estimate (\ref{eq:3.9}) holds for $\lambda$ belonging to the region
(\ref{eq:1.10}) with a suitable choice of the constant $C$. Thus, for such $\lambda$, we get $\|u_1\|=0$, 
which implies $u_1\equiv 0$, and hence $f\equiv 0$, as desired.

The proof of Theorem \ref{1.2} is much more complicated. We will need the next two lemmas. 
Let $\chi\in C^\infty(T^*\Gamma)$ be of compact support such that $1-\chi$ is supported in the region $\{r_0\ge\sigma\}$, 
where $\sigma\gg 1$ is a constant to be fixed in the next lemma. 

\begin{lemma} \label{7.2} 
For a suitable choice of $\sigma$ we have the estimate
\begin{equation}\label{eq:7.6}
\|{\rm Op}_h(1-\chi)f\|_{1,0}\lesssim h\|f\|_{0}
\end{equation}
for $\lambda\in\mathbb{C}^+\setminus\mathcal{L}_N^+(C_N)$ satisfying (\ref{eq:7.3}).
\end{lemma}

{\it Proof.} Choose a compactly supported function $\eta\in C^\infty(T^*\Gamma)$ such that $1-\eta$ is supported in the region 
$\{r_0\ge\sigma\}$
and $\chi=1$ on supp$\,\eta$. Define $\rho_j$, $\mathcal{H}_j$, $\mathcal{E}_j$, $j=1,2$, by replacing in the definition of $\rho$, 
$\mathcal{H}$, $\mathcal{E}$
in Sections 4 and 6 the functions $c,n$ by $c_j,n_j$. Clearly, taking $\sigma$ big enough we can arrange that the functions 
$1-\chi$ and $1-\eta$ are supported in both elliptic regions $\mathcal{E}_1$ and $\mathcal{E}_2$. 
Since $Tf=0$, we have the identity
$${\rm Op}_h((c_1\rho_1-c_2\rho_2)(1-\chi))f$$
$$+{\rm Op}_h(1-\chi){\rm Op}_h((c_1\rho_1-c_2\rho_2)(1-\eta))f-{\rm Op}_h((c_1\rho_1-c_2\rho_2)(1-\chi))f$$
$$={\rm Op}_h(1-\chi){\rm Op}_h((c_1\rho_1-c_2\rho_2)(1-\eta))f$$
$$={\rm Op}_h(1-\chi)\left(T{\rm Op}_h(1-\eta)+{\rm Op}_h((c_1\rho_1-c_2\rho_2)(1-\eta))\right)f$$
$$+{\rm Op}_h(1-\chi)T{\rm Op}_h(\eta)f.$$
Since $(c_1\rho_1-c_2\rho_2)(1-\chi)$ and $(c_1\rho_1-c_2\rho_2)(1-\eta)$ belong to $S^1(\Gamma)$, the $h-\Psi$DO calculus give
$${\rm Op}_h(1-\chi){\rm Op}_h((c_1\rho_1-c_2\rho_2)(1-\eta))-{\rm Op}_h((c_1\rho_1-c_2\rho_2)(1-\chi))
=O(h):H^1(\Gamma)\to L^2(\Gamma).$$
Therefore, using Theorem \ref{6.1} together with Lemma \ref{6.2} we obtain
\begin{equation}\label{eq:7.7}
\|{\rm Op}_h((c_1\rho_1-c_2\rho_2)(1-\chi))f\|_0\lesssim h\|f\|_{1,0}.
\end{equation}
On the other hand, we have
$$(c_1\rho_1+c_2\rho_2)(c_1\rho_1-c_2\rho_2)=c_1^2\rho_1^2-c_2^2\rho_2^2=(c_1^2-c_2^2)r_0-z(c_1n_1-c_2n_2).$$
Hence
$$\left|(c_1^2-c_2^2)r_0-z(c_1n_1-c_2n_2)\right|\lesssim\langle\xi'\rangle\left|c_1\rho_1-c_2\rho_2\right|.$$
On the other hand, in view of assumption (\ref{eq:1.2}) we have $|c_1^2-c_2^2|\ge C_0$ with some constant
$C_0>0$. Therefore, taking $\sigma$ big enough we can arrange that
$$\left|(c_1^2-c_2^2)r_0-z(c_1n_1-c_2n_2)\right|\ge C_1\langle\xi'\rangle^2,\quad C_1>0,$$
on supp$(1-\eta)$. Hence
\begin{equation}\label{eq:7.8}
\left|c_1\rho_1-c_2\rho_2\right|\ge C_2\langle\xi'\rangle,\quad C_2>0,
\end{equation}
on supp$(1-\eta)$, which implies
\begin{equation}\label{eq:7.9}
\|{\rm Op}_h(1-\chi)f\|_{1,0}\lesssim 
\|{\rm Op}_h((c_1\rho_1-c_2\rho_2)(1-\chi))f\|_0+h\|f\|_0.
\end{equation}
By (\ref{eq:7.7}) and (\ref{eq:7.9}).
$$\|{\rm Op}_h(1-\chi)f\|_{1,0}\lesssim h\|f\|_{1,0}\lesssim h\|f\|_{0}+h\|{\rm Op}_h(1-\chi)f\|_{1,0}.$$
We now absorb the last term in the right-hand side of the above estimate and get (\ref{eq:7.6}).
\eproof

By (\ref{eq:4.1}), (\ref{eq:7.5}) and (\ref{eq:7.6}) we obtain 
$$\|f\|_{1,0}\lesssim |\lambda|^{(\ell_2+1)/2}\tau_1(|\lambda|)^3\tau_2(|\lambda|)|{\rm Im}\,\lambda|^{1/2}\|f\|_{1,0}+h\|f\|_{0}.$$
Taking $h$ small enough we absorb the last term and arrive at the estimate
\begin{equation}\label{eq:7.10}
\|f\|_{1,0}\le C|\lambda|^{(\ell_2+1)/2}\tau_1(|\lambda|)^3\tau_2(|\lambda|)|{\rm Im}\,\lambda|^{1/2}\|f\|_{1,0},
\end{equation}
which holds for all $\lambda\in\mathbb{C}^+\setminus\mathcal{L}_N^+(C_N)$ satisfying (\ref{eq:7.3}) and (\ref{eq:7.4}), where
 $C>0$ is a constant independent of $\lambda$ and $N$. Hence, if 
\begin{equation}\label{eq:7.11}
\lambda\in\mathbb{C}^+\setminus\mathcal{L}_N^+(C_N),\,|{\rm Im}\,\lambda|\le (2C)^{-2}|\lambda|^{-(\ell_2+1)}\tau_1(|\lambda|)^{-6}\tau_2(|\lambda|)^{-2},
\end{equation}
 we can absorb the term in the right-hand side of
(\ref{eq:7.10}) and conclude that $\|f\|_{1,0}=0$, which implies $f\equiv 0$. In other words there are no
transmission eigenvalues in the region (\ref{eq:7.11}), as desired.

In what follows we will assume the conditions (\ref{eq:1.6}) and (\ref{eq:1.11}) fulfilled and we will show that in this case
the factor $|\lambda|^{1/2}$ in the right-hand side of
(\ref{eq:7.10}) can be removed.  Clearly, the condition (\ref{eq:1.6}) implies $\mathcal{H}_2\subset\mathcal{H}_1$
and $\mathcal{E}_1\subset\mathcal{E}_2$. Given a parameter $0<\varepsilon\ll 1$, independent of $\lambda$, we can choose
functions $\chi_\varepsilon^-, \chi_\varepsilon^0, \chi_\varepsilon^+\in C^\infty(T^*\Gamma)$ such that 
$\chi_\varepsilon^-+\chi_\varepsilon^0+\chi_\varepsilon^+\equiv 1$, supp$\,\chi_\varepsilon^-\subset\mathcal{H}_1$,
supp$\,\chi_\varepsilon^+\subset\mathcal{E}_1$, $\chi_\varepsilon^0=1$ in $\{|r_0/n_{\flat,0}-1|\le\varepsilon\}$ and 
$\chi_\varepsilon^0=0$ in $T^*\Gamma\setminus\{|r_0/n_{\flat,0}-1|\le 2\varepsilon\}$, where $n_{\flat,0}=\frac{n_1}{c_1}|_\Gamma$.
Clearly, supp$\,\chi_\varepsilon^+\subset\mathcal{E}_2$. 
Taking $\varepsilon$ small enough we can also arrange that supp$(1-\chi_\varepsilon^-)\subset\mathcal{E}_2$.
Using this we will prove the following

\begin{lemma} \label{7.3}
We have the estimate
\begin{equation}\label{eq:7.12}
\|{\rm Op}_h(1-\chi_\varepsilon^-)f\|_{1,0}\lesssim \left(|\lambda|^{\ell_2/2}\tau_1(|\lambda|)^5\tau_2(|\lambda|)|{\rm Im}\,\lambda|^{1/2}+h^{1/2}\tau_1(|\lambda|)^2+\varepsilon\right)\|f\|_{1,0}
\end{equation}
for $\lambda\in\mathbb{C}^+\setminus\mathcal{L}_N^+(C_N)$ satisfying (\ref{eq:7.3}) and (\ref{eq:7.4}).
\end{lemma}

{\it Proof.} Observe that the condition (\ref{eq:1.6}) implies the inequality
$$\frac{c_2n_2-c_1n_1}{c_2^2-c_1^2}<\frac{n_1}{c_1}.$$
Therefore the inequality (\ref{eq:7.8}) holds on $\mathcal{E}_1$. 
This implies the estimate (\ref{eq:7.9}) with $1-\chi$ replaced by $\chi_\varepsilon^+$.
On the other hand, since the function $\chi_\varepsilon^+$ is supported in both elliptic regions,
the estimate (\ref{eq:6.3}) holds with $1-\eta$ replaced by $\chi_\varepsilon^+$. 
Thus, in the same way as in the proof of Lemma \ref{7.2},
using (\ref{eq:5.4}) and (\ref{eq:7.5}) instead of (\ref{eq:6.15}), we conclude that 
\begin{equation}\label{eq:7.13}
\|{\rm Op}_h(\chi_\varepsilon^+)f\|_{1,0}\lesssim \left(|\lambda|^{\ell_2/2}\tau_1(|\lambda|)^5\tau_2(|\lambda|)|{\rm Im}\,\lambda|^{1/2}+h^{1/2}\tau_1(|\lambda|)^2\right)\|f\|_{1,0}.
\end{equation}
Choose a function  $\eta_\varepsilon\in C_0^\infty(T^*\Gamma)$
such that $\eta_\varepsilon=1$ in $\{|r_0/n_{\flat,0}-1|\le 3\varepsilon\}$ and 
$\eta_\varepsilon=0$ in $T^*\Gamma\setminus\{|r_0/n_{\flat,0}-1|\le 4\varepsilon\}$. Clearly, $\eta_\varepsilon=1$
on supp$\,\chi_\varepsilon^0$. Moreover, taking $\varepsilon$ small enough we can arrange that 
supp$\eta_\varepsilon\subset\mathcal{E}_2$. Since $Tf=0$, we have the identity
$${\rm Op}_h(\chi_\varepsilon^0\rho_2)f+\left({\rm Op}_h(\chi_\varepsilon^0){\rm Op}_h(\eta_\varepsilon\rho_2)f
-{\rm Op}_h(\chi_\varepsilon^0\rho_2)f\right)$$
$$={\rm Op}_h(\chi_\varepsilon^0){\rm Op}_h(\eta_\varepsilon\rho_2)f=-{\rm Op}_h(\chi_\varepsilon^0)\mathcal{N}_2(\lambda)
{\rm Op}_h(\eta_\varepsilon)f$$
$$+{\rm Op}_h(\chi_\varepsilon^0)\left(\mathcal{N}_2(\lambda){\rm Op}_h(\eta_\varepsilon)f
+{\rm Op}_h(\eta_\varepsilon\rho_2)f\right)$$
$$={\rm Op}_h(\chi_\varepsilon^0)\mathcal{N}_2(\lambda){\rm Op}_h(1-\eta_\varepsilon)f
+{\rm Op}_h(\chi_\varepsilon^0)c_2^{-1}c_1\mathcal{N}_1(\lambda,m)f$$
$$+{\rm Op}_h(\chi_\varepsilon^0)\left(\mathcal{N}_2(\lambda){\rm Op}_h(\eta_\varepsilon)f
+{\rm Op}_h(\eta_\varepsilon\rho_2)f\right).$$
The $h-\Psi$DO calculus give
$${\rm Op}_h(\chi_\varepsilon^0){\rm Op}_h(\eta_\varepsilon\rho_2)-{\rm Op}_h(\chi_\varepsilon^0\rho_2)
=O_\varepsilon(h):L^2(\Gamma)\to L^2(\Gamma),$$
$${\rm Op}_h(\chi_\varepsilon^0)c_2^{-1}c_1-c_2^{-1}c_1{\rm Op}_h(\chi_\varepsilon^0)
=O_\varepsilon(h):L^2(\Gamma)\to L^2(\Gamma).$$
Therefore, using the estimates (\ref{eq:4.11}), (\ref{eq:5.1}), (\ref{eq:5.4}), (\ref{eq:6.3}), (\ref{eq:6.14}) (with 
$m\equiv 0$) and (\ref{eq:7.5}), we get
\begin{equation}\label{eq:7.14}
\|{\rm Op}_h(\chi_\varepsilon^0\rho_2)f\|_0\lesssim \left(|\lambda|^{\ell_2/2}\tau_1(|\lambda|)^5\tau_2(|\lambda|)|{\rm Im}\,\lambda|^{1/2}+h^{1/2}\tau_1(|\lambda|)^2+\varepsilon\right)\|f\|_{1,0}.
\end{equation}
On the other hand, the condition (\ref{eq:1.6}) guarantees that $|\rho_2|\ge C>0$ on supp$\,\chi_\varepsilon^0$,
provided $\varepsilon$ is taken small enough. Hence
\begin{equation}\label{eq:7.15}
\|{\rm Op}_h(\chi_\varepsilon^0)f\|_0\lesssim 
\|{\rm Op}_h(\chi_\varepsilon^0\rho_2)f\|_0+h\|f\|_0.
\end{equation}
Clearly, the estimate (\ref{eq:7.12}) follows from (\ref{eq:7.13}), (\ref{eq:7.14}) and (\ref{eq:7.15}).
\eproof

By (\ref{eq:4.5}), (\ref{eq:7.5}) and (\ref{eq:7.12}) we obtain 
$$\|f\|_{1,0}\lesssim \left(|\lambda|^{\ell_2/2}\tau_1(|\lambda|)^5\tau_2(|\lambda|)|{\rm Im}\,\lambda|^{1/2}+|\lambda|^{-1/2}\tau_1(|\lambda|)^2+\varepsilon\right)\|f\|_{1,0}.$$
Taking $|\lambda|^{-1}$ and $\varepsilon$ small enough, in view of the assumption (\ref{eq:1.11}), we can absorb the last two terms to obtain
\begin{equation}\label{eq:7.16}
\|f\|_{1,0}\le C|\lambda|^{\ell_2/2}\tau_1(|\lambda|)^5\tau_2(|\lambda|)|{\rm Im}\,\lambda|^{1/2}\|f\|_{1,0}
\end{equation}
which holds for all $\lambda\in\mathbb{C}^+\setminus\mathcal{L}_N^+(C_N)$ satisfying (\ref{eq:7.3}) and (\ref{eq:7.4}), where
 $C>0$ is a constant independent of $\lambda$ and $N$. Hence, if 
\begin{equation}\label{eq:7.17}
\lambda\in\mathbb{C}^+\setminus\mathcal{L}_N^+(C_N),\,|{\rm Im}\,\lambda|\le (2C)^{-2}|\lambda|^{-\ell_2}
\tau_1(|\lambda|)^{-10}\tau_2(|\lambda|)^{-2},
\end{equation}
 we can absorb the term in the right-hand side of
(\ref{eq:7.16}) and conclude that $\|f\|_{1,0}=0$, which implies $f\equiv 0$. In other words there are no
transmission eigenvalues in the region (\ref{eq:7.17}) in this case, which is the desired conclusion.

\end{document}